\magnification=\magstep1
\baselineskip=16pt
\lineskip=0pt
\lineskiplimit=0pt
\parskip=6pt
\vsize=7.65in
\font\small = cmr10 at 6pt

\input epsf.tex

% Special Roman or Blackboard Bold

\def\Z   {{\bf Z}}

\def\R   {{\bf R}}
\def\C   {{\bf C}}
\def\P   {{\bf P}}

\def\dist {\atop {\hbox{\small distinct}}}

% Roman

\def\Re    {{\rm Re}}
\def\Im    {{\rm Im}}
\def\tr    {{\rm tr}}

% For this paper

\def\Ktilde    {\tilde{\cal K}}

\def\Kahler    {K\"ahler\ }
\def\CY        {Calabi--Yau\ }

\def\d         {\partial}
\def\dbar      {\bar\partial}
\def\pi        {\partial /\partial z_i}
\def\pj        {\partial /\partial z_j}

\def\cp				    {\scriptstyle \cup \displaystyle}
\def\ip {\, \raise 1.5pt \hbox{\underbar{\hskip 4.3pt }\vrule height 1.5pt
depth 1.4pt width .2pt}\, }

\def\ha        {\textstyle {1\over 2}\displaystyle}
\def\hai        {\textstyle {i\over 2}\displaystyle}
\def\th        {\textstyle {1\over 3}\displaystyle}
\def\qu        {\textstyle {1\over 4}\displaystyle}
\def\si        {\textstyle {1\over 6}\displaystyle}
\def\iha       {\textstyle {i\over 2}\displaystyle}
\def\tth       {\textstyle {2\over 3}\displaystyle}

% Caligraphic Letters

\def\calH      {{\cal H}}

\def\calK      {{\cal K}}

\def\calH      {{\cal H}}

\noindent \bf Sectional curvatures of \Kahler moduli\rm
\bigskip
\bigskip  \noindent \rm P.M.H. Wilson \medskip
\noindent Department of Pure Mathematics, University of Cambridge,\hfil\break 16
Wilberforce Road, Cambridge CB3 0WB, UK
\hfil\break email : pmhw@dpmms.cam.ac.uk \medskip
\noindent Max-Planck-Institut f\"ur Mathematik,\hfil\break
 Vivatsgasse 7, 53111 Bonn, Germany

\bigskip

\noindent \bf Abstract 

\rm  If $X$ is a compact \Kahler manifold of dimension $n$,
we let $\calK \subset H^{1,1} (X, \R )$ denote the cone of 
\Kahler classes, and $\calK _1$ the level set given 
by classes $D$ with $D ^n = 1$.  This space is 
naturally a Riemannian manifold 
and is isometric to the manifold $\tilde \calK _1$ of 
\Kahler \it forms \rm $\omega$ with $\omega ^n$ some fixed volume form, equipped with the Hodge metric, as studied previously by Huybrechts. 
We study these spaces further, in particular their geodesics and sectional curvatures.  Conjecturally, at least for \CY 
manifolds and probably rather more generally, these sectional curvatures should be  bounded between $-\ha n(n-1)$ and zero.   We find  simple formulae 
for the sectional curvatures, and prove both the bounds hold 
for various classes of varieties, 
developing along the way a mirror to 
the Weil--Petersson theory of complex moduli.  In the case of
threefolds with $h^{1,1} =3$, we produce an explicit formula for this curvature in terms of the invariants of the cubic form.  This enables us to check the bounds by computer for a wide range of examples.    Finally, we  
explore the implications of the non-positivity of these curvatures.

\bigskip

\noindent \bf{0.\quad Introduction} \rm

\medskip \rm  Let $X$ denote a compact \Kahler manifold of dimension $n$.
The cup product on $H^2 (X, \Z )$ determines a degree $n$ form  
on $H^{1,1} (X, \R )$, and we define the
 \it positive cone \rm to be  
$ \{ D \in H^{1,1} (X, \R ) \, :\, D^n >0 \} $.   Cup product also determines
an \it index cone \rm (cf. [18]), which will be denoted by
$ W$, consisting of elements $D$ in the positive cone for which the quadratic form 
on $H^{1,1} (X, \R )$ given by $L \mapsto D^{n-2}\cp L^2 $ has signature $(1, h^{1,1} -1 )$.  By
the Hodge index theorem, we have an inclusion of open cones 
in $H^{1,1} (X, \R )$ from the \Kahler cone 
$\calK$ into $ W$.  We 
denote by $W_1$ the level set $\{ D\in W \, : \, D^n = 1 \}$ in the index cone.  It is 
clear that $W_1$ is a smooth manifold, usually non-connected,
whose tangent space at a point $D$ is
identified as $\{ L \in H^{1,1} (X, \R ) \, : \, D^{n-1}\cp L =0 \}$.  Moreover, there is a
natural Riemannian metric on $W_1$ given by the pairing, for tangent vectors $L_1 , L_2$ 
at $D \in W_1$,
$$ (L_1 , L_2) \mapsto - D^{n-2}\cp L_1 \cp L_2 .$$
We note that for $n=2$, the cup product on $H^{1,1} (X, \R)$ is just a Lorentzian real quadratic
form, and $W$ coincides with the positive cone.
Moreover, on each connected component of $W_1$,  our construction
reduces to the standard construction 
of real hyperbolic space 
(Example 10.2 in Chapter XI of [11], or page 189 of [8]); in particular it has constant negative curvature $-1$.

For a given complex structure on $X$, the \Kahler structure is determined by 
the \Kahler form $\omega$, a closed real
$(1,1)$-form.  
For a fixed complex structure and volume form $ \omega _0 ^n /n! $, Huybrechts
introduced what he called the \it curved 
\Kahler cone \rm $\Ktilde$ consisting of all \Kahler forms $\omega$ with $\omega ^n = c\,
\omega _0 ^n$ for some $c>0$ [7].  The Aubin--Calabi--Yau theorem [1] then implies
that  the projection map from $\Ktilde$ to $\calK $ is a bijection.
We normalise the volume form so that $\int _X 
\omega _0 ^n  = 1$, and set $\tilde \calK _1$ to consist of the 
\Kahler forms $\omega$ with $\omega ^n = 
\omega _0 ^n$. Thus,  setting $\calK _1 = \calK \cap
W_1$,
the projection $\tilde \calK _1 \to \calK _1 $ is a bijection.  Moreover,  
Huybrechts observes that $\Ktilde _1$ is a smooth manifold, with the tangent space 
at $\omega \in \Ktilde _1$ consisting of the \it primitive \rm closed real $(1,1)$-forms
$\alpha$;
 essentially this is saying that to first order $(\omega + \epsilon\, \alpha )^n = \omega
^n$. He remarks that the Hodge identities imply that such forms are harmonic with 
respect to $\omega$, and so the tangent space is identified as the 
$(h^{1,1} - 1)$-dimensional space of primitive harmonic $(1,1)$-forms.  There is now a
natural Riemannian metric on $\Ktilde _1 $ given by the Hodge metric on harmonic forms.
Given elements $\alpha _1 , \alpha _2$ in the tangent space to $\Ktilde _1$ at $\omega$,
the metric is specified by the pairing 
$$(\alpha _1 , \alpha _2) \mapsto - \int _X \omega ^{n-2}\wedge \alpha _1 \wedge \alpha
_2.$$  Thus the projection map $\Ktilde _1 \to \calK _1 $ is in fact an isometry of 
Riemannian manifolds, enabling us to identify $\Ktilde _1$ with $\calK _1 $. 
We shall call this Riemannian manifold 
the \it normalised \Kahler moduli space\rm , and it is now clearly independent of the
choice of normalised volume form.

In the case when $X$ is a \CY $n$-fold, not necessarily with $h^{2,0} =0$,
we have a nowhere vanishing holomorphic 
$n$-form $\Omega$ on $X$, and we can take the volume form to be a suitable multiple
of $(i/2)^n \,\Omega \wedge \bar\Omega $.  In this case, one 
is struck by the similar properties enjoyed by $\calK _1$ and the space of complex
structures on a \CY $n$-fold.  For instance, $\calK _1$ has degenerations at \it finite
distance\rm , which correspond to singular \CY $n$-folds with canonical singularities 
(cf [27] in dimension 3), and also degenerations at \it infinite distance.  \rm This may
be  compared with the properties of the Weil--Petersson metric for the complex moduli 
space [25].  In Section 1, we shall use Mirror Symmetry to suggest an
explanation  for such similarities.  The sectional curvatures of the
Weil--Petersson metric on the complex moduli space of the mirror are expected to be non-positive near the large complex structure limit;
 we argue why this might indicate  that the sectional curvatures
on $\calK _1$ are also non-positive.  In general, we shall 
say that the  \it \Kahler moduli curvature 
is semi-negative 
\rm if the sectional curvatures of $\calK _1$ are non-positive.

The argument from Section 1 suggests that a mirror version of Weil--Petersson theory
should be developed, to hold on the normalised \Kahler moduli space; 
the basics of such a
theory are developed in Sections 2 to 4.  In Section 2, a   bracket 
operation $A^{1,1} \times A^{1,1} \to A^{2,1}$ is defined, a 
mirror version of the Tian--Todorov Lemma is
proved, and the structure of a differential graded Lie algebra is defined 
on $\bigoplus _{p\ge 0} A^{p,1}$ 
(throughout the paper, $A^{p,q} (X)$ denotes the space of $(p,q)$-forms on $X$).
  In Section 3, we investigate geodesics on the  normalised \Kahler moduli space, finding the all important quadratic term (3.6), and
in Section 4, we obtain a simple formula (4.1) for the sectional curvatures, 
which is relevant for the question of semi-negativity, and another formula (4.3), which leads us naturally to the
question of a \it lower \rm bound.  The following question 
turns out to be the natural one:

\proclaim Question.  For which compact \Kahler manifolds
$X$ are the sectional curvatures of $\calK _1$ non-positive, and for
which $X$ are they bounded below by 
 $-\ha n(n-1)$?
\par

We conjecture that both these bounds should hold
for \CY manifolds, but 
in fact they likely to be true rather more widely 
than that.    The author 
knows of no examples of compact \Kahler manifolds where the bounds fail.
The question of when the formulae (4.1) and (4.3) 
do yield the suspected bounds, and also a strategy for proving them, are discussed in (4.4).

An instructive example where the bounds hold 
 is provided (4.2) when $X$ is a complex torus.  We observe that 
the normalised \Kahler moduli space
$\calK _1$ may be identified as the space of positive definite
hermitian matrices of determinant 1, that is, the symmetric space $ SL(n, \C )/SU(n)$, 
and that the Hodge metric on $\calK _1$ corresponds to some multiple of the symmetric space metric.
The standard formula for the sectional curvatures of $ SL(n, \C )/SU(n)$ then 
corresponds to the formula we derive in (4.1).

If we could prove the conjectured semi-negativity  
 for \CY manifolds, it would provide new information on
the possible location of the \Kahler cone $\calK$ in cohomology, and potentially useful  
information  concerning which differentiable manifolds may support \CY structures. 
  For $h^{1,1} (X)
\le 2$, the conjecture gives no information, but one should expect increasingly more
information as $h^{1,1}$ increases.    

The case of \Kahler manifolds with 
$h^{1,1} = 3$ does have the advantage that the curvature is easier to calculate, and for
threefolds we derive a very explicit formula (5.1) for it in terms of the invariants
of the ternary cubic form given by cup product (this formula has recently been shown to have a natural extension, 
in a similar shape, for arbitrary degrees $>2$ [23]).
Using the formula, the author has checked a large number of \Kahler threefolds  with $b_2 = 3$, and  has verified the bounds in each case. 
  In the case for instance of complete 
intersections in the product of three projective spaces, there is persuasive 
numerical evidence for the bounds; here, the lower bound 
should be $-9/4$ rather than $-3$, and this fact has been
verified by computer for the case of complete intersections
in $\P ^5 \times \P ^5 \times \P ^5$.  The upper bound 
has been verified for the case of complete 
intersections in $\P ^3 \times \P ^2 \times \P ^2$.

In these examples,  
the semi-negativity condition usually places stronger restrictions concerning the 
location of the \Kahler cone than are provided by just the index cone.  
On the other hand, we show 
 that,  for \Kahler threefolds with $b_2 (X) =3$, even the stronger condition 
that the \Kahler moduli curvature lies between $-3$ and zero 
does not rule out further cases for the cubic form than are excluded by 
the standard results from K\"ahler geometry.  I conjectured in an earlier version
of this paper that this would not be the case for higher dimensions or higher 
second Betti number, and this expectation has recently been confirmed in [23].  
  The case of threefolds is of particular interest
because, by results of Wall and Jupp, there is a very simple criterion for a given 
cubic form to be realisable as the cup product on $H^2 (X, \Z )$, for some smooth 
simply connected 6-manifold $X$ with torsion free homology, and moreover such
manifolds are classified by means of their invariants [18].  Note that all simply 
connected 6-manifolds are formal [16], and so in dimension 6 we cannot use this to
distinguish those manifolds not supporting a \Kahler structure. 
If the conjectured semi-negativity is true, then we can 
 give examples of simply connected compact differentiable 6-manifolds 
with torsion free homology (specified by invariants
including the cubic form, the first Pontryagin class and the third betti number) which do not
support any \CY structures.  The non-existence of \CY structures 
on these 6-manifolds would not appear to follow from existing criteria.\medskip 

Finally, I mention two cases where the bounds on  
the \Kahler moduli sectional curvatures are easy to check.

\noindent\it Example 1.\quad\rm  Suppose $X$ is a compact irreducible complex
symplectic (hence hyperk\"ahler) manifold of dimension $2n$.
Letting $q_X$ denote the Beauville--Bogomolov 
quadratic form, this has the property that for any 
$\alpha \in H^2 (X, \R )$, $q_X (\alpha ) ^n = \lambda 
\int _X \alpha ^n $, for some fixed positive constant $\lambda$.  Moreover, 
$q_X$ is positive on any \Kahler class $\omega$ 
and negative definite on the corresponding primitive classes
$H^{1,1} (X, \R )_\omega $.  Setting $2\, b(\alpha ,\beta ) = q_X (\alpha 
+ \beta ) - q_X (\alpha ) - q _X (\beta )$, we observe that for
$\omega$ a \Kahler class and $\alpha$ a primitive $(1,1)$-class, 
$q _X (\omega + t \alpha ) = q_X (\omega ) + 2t\, b(\omega ,\alpha) + t^2 q_X (\alpha ),$
whilst $\int _X (\omega + t\, \alpha )^n = \int _X \omega ^n + n(2n-1) t^2\, 
\omega ^{2n-2}\cp \alpha ^2 + O(t^3 );$  thus $b(\omega , \alpha )=0$.
 Assuming that $\int \omega ^n =1$ and that $q_X$ is normalised so that $\lambda =1$, we have $q_X(\omega )=1$, and we deduce that 
$q_X (\alpha ) = (2n-1)\, \omega ^{2n-2} \cp 
\alpha ^2 $ for all $\alpha \in H^{1,1}(X, \R )_\omega $.  Since $q_X$ is Lorentzian on $H^{1,1}(X, \R )$, this implies that $W_1$ is hyperbolic, as in the surface case,
with constant sectional curvatures of value 
$-(2n-1)$.  Note that this is consistent with the fact that, for a general irreducible complex symplectic manifold, the \Kahler cone is a connected component of what we've called the positive cone.
\medskip
\noindent\it Example 2.\quad\rm  Suppose $h^{1,1} (X) = m+1$, and there are $m$
independent divisorial contractions of irreducible divisors $E_i$ on $X$ to (distinct) points. 
We have a basis of $H^{1,1} (X, \R ) \subset H^2 (X, \R )$ 
given by $H, -s_1 E_1 , \ldots , - s_m E_m$, with $H$ nef and the $s_i >0$, 
with respect to which 
the degree $n$ form on $H^{1,1} (X, \R )$ is diagonal of the form 
$$ x_0 ^n - x_1 ^n - \ldots - x_m ^n .$$  
Let $W_1 ^+$ denote the open subset of $W_1$ 
given by intersecting with the
positive orthant (with all coordinates strictly positive); when $n>2$, this is just the
connected component of $W_1$ that contains $\calK _1$.  The metric on $W_1 ^+$ is the
restriction of  
$$ - x_0 ^{n-2} dx_0 ^2 + x_1 ^{n-2} dx_1 ^2 + \ldots + x_m ^{n-2} dx_m ^2 .$$

Let $U_1$ denote the upper sheet of the hypersurface $ y_0 ^2 - y_1 ^2 - \ldots - y_m ^2 =1$ 
in $\R ^{m+1}$,
equipped with the metric given by restricting  
$ -  dy_0 ^2 +  dy_1 ^2 + \ldots +  dy_m ^2 $; as commented above, this is 
a standard construction for 
$m$-dimensional hyperbolic space.  We let $U_1 ^+$ denote the intersection of $U_1$ with the
positive orthant, and consider the diffeomorphism from $W_1 ^+$ onto $U_1 ^+$ given by 
$y_i = x_i ^{n/2}$.  Since $dy_i = {n \over 2}\, x_i ^{(n-2)/2}\, dx_i$, the pullback of the
hyperbolic  metric to $W_1 ^+$ is simply the given metric scaled by a factor of $({n\over
2})^2$.  So scaling distances on $W_1 ^+$ by a factor of ${2\over n}$, we obtain a manifold
isometric to 
$U_1 ^+$.  Since hyperbolic space has constant sectional curvatures of value $-1$, our manifold
$W_1 ^+$ has constant sectional curvatures of value $-({n\over 2})^2$.  For $n=3$, this is the
reason for the term $-9/4$ in  (5.1).

Some of the research on this paper was pursued by the author during a visit to the
University of California at San Diego, as a guest of Mark Gross, whom he would like to
thank both for his hospitality and useful discussions.  He would like to thank 
UCSD for their financial support, and Burt Totaro for his interest and 
insightful comments on a preliminary version of the manuscript.  He also thanks 
Mike Douglas and Zhiqin Lu for comments concerning the curvature of the 
Weil--Petersson metric.  The paper was
completed during a stay at the Max-Planck-Institut in Bonn, 
which he would like to thank 
for financial support.

\bigskip

\noindent \bf{1.\quad Motivation from Mirror Symmetry} \medskip \rm

\medskip \rm One possible interpretation of Mirror Symmetry is in terms of the SYZ 
Conjecture.  The simplest case here is that of a (non-compact) \CY $n$-fold $X$ fibered by 
special lagrangian tori for which the Ricci flat metric is semiflat (that is invariant
under  the torus action).  This case was studied by Hitchin, Gross and Leung [6,5,12].  
If the holomorphic 
$n$-form on $X$ is denoted by $\Omega$ and the \Kahler form by $\omega$, recall that the 
condition of Ricci flatness may be written as
$$ \omega ^n /n! = c\, (\hai )^n \, \Omega \wedge \bar\Omega $$
for some positive real constant $c$.  The solutions sought correspond to solutions (with
appropriate boundary conditions) of the real Monge--Amp\`ere equation on the base of the
fibration.

In [12], 
Leung produces an explicit holomorphic map (defined in terms of the fibration) between the 
\Kahler moduli space of $X$ and the complex moduli space of the mirror.  In Leung's
notation, fixing a standard complex structure on $X$, the \Kahler moduli space
consists of 
$T^n$-invariant complex
2-forms $\omega + i\beta$, with $\omega$ as above, and furthermore satisfying 
$$ \Im \, (e^{i\theta} (\omega + i \beta )^n ) = 0$$
for some angle $\theta$.  At a given point $\omega + i\beta$ of this space, Leung remarks
that the  tangent space consists of complex harmonic $(1,1)$-forms, and 
he argues that therefore
the  Hodge $L^2$-metric on the harmonic forms (defined at a given point by the real part
$\omega$) determines a Riemannian metric on the whole 
space.  On the corresponding space of complex structures on the mirror, one can take the
Weil--Petersson metric  (also defined as a Hodge $L^2$-metric on harmonic forms), and he
argues formally that the explicit correspondence that he has defined is an isometry (see
[12], Section 2).

For compact \CY $n$-folds, we cannot hope for the metrics to be semi-flat and there will be
instanton corrections to consider.  We can however argue as follows: 
by standard theory, a \CY manifold in the sense used in this paper has
a finite unramified cover which is the 
product of a complex torus, irreducible
complex symplectic manifolds, and 
simply connected \CY manifolds with 
$h^{2,0} = 0$.  As we know the 
conjectured bounds hold for the 
first two types of manifold, we may reduce the question to
the third type; we do not prove this assertion here, but it follows
for instance from the formula derived in (4.1), the key point here
being that the \it mixed \rm sectional curvatures of $\tilde\calK _1$
for the product are necessarily zero.  
For simply connected \CY manifolds in dimension 3 with $h^{2,0} = 0$, the 
Main Theorem from [27] says that 
 the \Kahler cone $\calK$ is essentially
invariant under deformations of the complex structure; 
to be precise, there 
are some complex codimension one loci in moduli, corresponding to
the existence of elliptic quasi-ruled surfaces, and here 
the \Kahler cone jumps down.  A similar result is true for general
$n$, without the explicit description of the loci where the 
\Kahler cone jumps down. 
The normalised \Kahler moduli
space $\calK_1$, considered as a subset of $H^2 (X, \R)$, 
can then be considered as an invariant of the complex 
structure.  When addressing the conjecture on 
the sectional curvatures of $\calK_1$, we are 
therefore at liberty to take 
arbitrary general points in the complex
moduli space, 
for instance degenerating to a large complex structure limit 
point on the boundary, if such a limit exists; 
 here the SYZ Conjecture will be relevant.
One hopes that Leung's results are then indicative of
what we might expect in the large
\Kahler structure limit.  Todorov [22] claims that on 
the complex moduli space of \CY $n$-folds, the
Weil--Petersson metric has non-positive sectional curvatures, but unfortunately the calculations on page 65 of [2] show this not to be true.  Nevertheless, it still seems likely that this is true near a large complex structure limit point.  For 
one dimensional moduli, it is shown in [26] that the Weil--Petersson metric 
is exponentially asymptotic to a scaling of the Poincar\'e metric for any 
degeneration at infinite distance, and in 
particular has negative curvature there.  The 
arguments only involve the variation of Hodge structure, and so it is 
reasonable to believe that that a similar statement is true near a 
large complex structure limit point (with the usual conventions of 
not getting too close to the complex codimension one 
discriminant loci meeting there) in 
the case of higher dimensional moduli.  The 
asymptotic behavior of the Weil--Petersson metric in the case of higher 
dimensional moduli is rather subtle, and is the subject of active and ongoing 
research (see for instance [14,15,3]).  
One is however led to conjecture from this that on the mirror side, 
 the Hodge metric on the complex \Kahler moduli space will have non-positive 
sectional curvatures near the 
large \Kahler structure limit.  On Leung's complex 
\Kahler  moduli space, complex conjugation defines an isometry.  The moduli space of
real
\Kahler forms is then just the fixed locus of this isometry, and is therefore a totally
geodesic submanifold (see [10], page 59).  Thus the real \Kahler moduli space should also 
have non-positive sectional curvatures, at least in the large \Kahler structure limit.

Suppose now $\omega \in \calK _1$, and let $g$ denote the Hodge metric on the tangent space
to $\calK _1$ at $\omega$, namely the primitive real harmonic $(1,1)$-forms.  Clearly then,
for $r>0$, 
the metric at $r \omega$  is $r( dr^2 + g)$.  For $\lambda >0$, there is a scaling 
map $\theta _\lambda $ from $\calK = \R_+ \times \calK _1$ to itself, given by 
$\theta _\lambda (r, x) = (r', x)$, where $r' = \lambda r$.  Moreover, we check that 
$$ \lambda ^{-1} \theta _\lambda ^* \, r'((dr') ^2 + g) = r(\lambda ^2 dr^2 + g).$$
Thus the scaling map sends the level set $r=1$ to the level set $r' = \lambda$, and the
normalised pullback of the metric from points on the level set $r' = \lambda$ is the 
same as the metric at the level set $r=1$, except that the distances in the radial direction
have been stretched by a factor $\lambda$.  Hence, in the limit, only the 
sectional curvatures along
the level set $\calK _1$ will survive.  In particular, the Hodge metric on $\calK _1$
should also have non-positive sectional curvatures.

All this is of course only a plausability argument. To prove that the \Kahler moduli curvature is semi-negative in 
the \CY case, the author anticipates that a theory 
describing 
\Kahler moduli which is mirror to the Weil--Petersson theory on complex moduli will be needed,
and the basics of such a theory are developed in the next three sections.

\bigskip

\noindent \bf{2.\quad Bracket operation on $(1,1)$-forms}
\medskip

\rm Let $X$ denote a compact \Kahler $n$-fold, equipped with a fixed 
complex structure  and a \Kahler structure, determined by a closed real
$(1,1)$-form $\omega$.  We shall assume that $\omega$ is normalised so that
$\int _X \omega ^n  =1$.

Motivated by the theory on the mirror side [21,22], 
we define a bracket 
operation on the $(1,1)$-forms $A^{1,1}$.  Given $\alpha \in A^{1,1}$, there exists
a unique element $\theta \in A^{1,0} (\Theta )$ such that $\theta \ip \omega =
\alpha$, where $\Theta$ denotes the sheaf of holomorphic vector fields on $X$ and $\ip$
denotes interior product.  At any given point, we can find 
local coordinates $z_1, \ldots
,z_n$,   for which $\omega = {i\over 2} \sum d z_i \wedge d \bar z_i$ and 
$\alpha = {i\over 2} \sum a_{ij} d z_i \wedge d \bar z_j$, and then $\theta$ can be
written at that point as  
$$\theta = \sum a_{ij}\ \pj \otimes dz_i.$$  Note that $\alpha$ is real if and only if 
$a_{ji} = \bar a_{ij}$ for all $i,j$ --- mostly we shall only consider real
$(1,1)$-forms.

If $\alpha$, $\beta$ are real $(1,1)$-forms, then with coordinates as above such that 
$\alpha = {i\over 2} \sum a_{ij}\ d z_i \wedge d \bar z_j$ and 
$\beta = {i\over 2} \sum b_{kl}\ d z_k \wedge d \bar z_l$,  we have at
the given point that 
$$ \theta \ip \beta =  {i\over 2} \sum a_{ij} b_{jl}\ d z_i \wedge d \bar z_l.$$
We observe for future use the identity 
$ \langle \alpha , \beta \rangle  = \Lambda (\theta \ip \beta ),$
where $\langle \, ,\, \rangle $ denotes the induced (pointwise) metric on forms and 
$\Lambda$ denotes the adjoint of the Lefshetz 
operator $L$, given in coordinates on page 114 of [4].
We also note that if $\beta = \phi \ip \omega $, then 
$$ \overline {\theta \ip \beta} = 
\bar\theta \ip \beta = \phi \ip \alpha ,$$
where $\bar \theta$ is an element of $A^{0,1} (\bar \Theta )$.

We can define a
(super) bracket operation
$$ [\ ,\ ]\quad : \quad A^{1,0} (\Theta ) \times A^{1,0} (\Theta ) \to 
A^{2,0} (\Theta),$$ or equivalently, 
via $[\alpha , \beta ] := [\theta , \phi ]\ip \omega $,
 a (symmetric) bracket 
$$ [\ ,\ ]\quad : A^{1,1} \times A^{1,1} \to A^{2,1}.$$  
The \Kahler condition enables us, at a given point, 
to choose 
local coordinates $z_1, \ldots ,z_n$ such that 
$$\omega = {i\over 2} \sum d z_i \wedge d \bar z_i \ + O(\, |z|^2).$$
We write $[\alpha , \beta ]$ at that
point as
$$ [\alpha , \beta ] = {i\over 2} \sum _{i,j,k,l} \bigg (a_{ij} {{\partial
b_{kl}}\over  {\partial z_j}} + b_{ij} {{\partial a_{kl}}\over {\partial z_j}}\bigg ) \
d \bar z_l
\wedge  d z_i \wedge d z_k ;$$ the corresponding 
bilinear form on $ A^{1,0} (\Theta )$ wedges the $dz_i$
and takes the Lie bracket of vector fields.
Since only one derivative is
being taken, this is independent of our choice of local 
coordinates of the specified type, a fact also confirmed by
(2.2) below, which provides a visibly coordinate-independent
characterization of the pairing.

\proclaim Lemma 2.1.  For $\alpha$ any
 $(1,1)$-form, $[\omega , \alpha ] = \partial\alpha$.\par

\noindent \it Proof.\rm\quad  Immediate from the local formula given above.\medskip

\proclaim Proposition 2.2.
With the notation as above,
$$ [\alpha , \beta ] \, + \,\partial (\theta \ip\beta\, + \,\phi \ip\alpha) 
\, = \, (\# \partial \alpha)\ip\beta \, + \, (\# \partial \beta)\ip\alpha \,
+ \, \theta\ip\partial\beta \, + \, \phi\ip\partial\alpha ,$$
where $\# : A^{2,1} \to A^{2,0}(\Theta )$ is the natural map given by 
$(\#\Gamma )\ip\omega = \Gamma$.\par

\noindent\it Proof.\quad\rm This 
 is the mirror of what is sometimes
referred to as the Tian--Todorov Lemma, and as in that case, once one knows the correct
statement of the result, the proof reduces to a (rather unenlightening) calculation in
terms of local coordinates.  If we wish to check the identity at a given point, we
choose local coordinates $z_1, \ldots ,z_n$ such that \break $\omega = {i\over 2}\sum
dz_i \wedge d\bar z_i
\ + O(\, |z|^2)$.  
At the given point therefore 
$$ [\alpha , \beta ] = {i\over 2} \sum _{i,j,k,l} \bigg (a_{ij} {{\partial
b_{kl}}\over  {\partial z_j}} + b_{ij} {{\partial a_{kl}}\over {\partial z_j}}\bigg ) \
d \bar z_l
\wedge  d z_i \wedge d z_k $$
and
$$ \partial (\theta \ip\beta ) =
{i\over 2} \sum _{i,j,k,l} \bigg (a_{ij} {{\partial
b_{jl}}\over  {\partial z_k}} + b_{jl} {{\partial a_{ij}}\over {\partial z_k}}\bigg ) 
\ d z_k
\wedge  d z_i \wedge d \bar z_l . $$

Also, since $\partial \alpha = {i\over 2} \sum _{i,j,k} {{\partial a_{ij}}\over {\partial z_k}} 
\ d z_k \wedge  d z_i \wedge d \bar z_j$, we see that 
$\# \partial \alpha = \sum _{i,j,k} {{\partial a_{ij}}\over {\partial z_k}} 
\ \pj \otimes d z_k \wedge  d z_i$ at the given point, 
from which it follows that
$$  (\# \partial \alpha )\ip \beta \ = \ {i\over 2}
\sum _{i,j,k} {{\partial a_{ij}}\over {\partial z_k}} b_{jl} 
\ d z_k \wedge  d z_i \wedge d \bar z_l . $$
Finally we have
$$ \theta \ip \partial\beta \ = \ 
{i\over 2} \sum _{i,j,k,l} \bigg (a_{ij} {{\partial
b_{kl}}\over {\partial z_j}} + a_{ij} {{\partial b_{jl}}\over {\partial z_k}}\bigg ) 
\ d z_i
\wedge  d z_k \wedge d \bar z_l . $$

From this it follows that, at the given point, 
$$ - \partial (\theta \ip\beta ) 
+ (\# \partial \alpha)\ip\beta \  
+ \ \theta\ip\partial\beta  \ = \ 
{i\over 2} \sum _{i,j,k,l} a_{ij} {{\partial
b_{kl}}\over  {\partial z_j}} \ 
d  z_i
\wedge  d z_k \wedge d \bar z_l .$$
The symmetry between $\alpha$ and $\beta$ now yields the formula claimed.

\proclaim Corollary 2.3.  If $\partial\alpha = 0 = \partial\beta$, then 
$ [\alpha , \beta ] + \partial (\theta \ip\beta\, + \, \phi \ip\alpha) =0$.

When $\alpha$, $\beta$ are real, 
the conditions of (2.3) are equivalent to $d\alpha = 0 = d\beta $.
 Note that by a previous comment, $\theta \ip\beta\, + \, \phi \ip\alpha\, = 2\, \Re (\theta \ip\beta)$ is 
then also a real $(1,1)$-form.
\bigskip

The pairing we've defined may be extended in an obvious way to give 
a bilinear pairing $$ A^{p,1} \times A^{q,1} \to A^{p+q, 1},$$ which is 
symmetric or antisymmetric according to the parity of $pq +1$, where the
corresponding pairing on $A^{p,0}(\Theta ) \times A^{q,0}(\Theta )$ is again defined by wedging the $dz_i$ and taking the Lie bracket of vector fields.  We shall also 
assume here that the metric is real analytic.

\proclaim Theorem 2.4.  The bracket defined above satisfies the Jacobi
identity
$$ [\gamma , [ \alpha , \beta ]] + [\beta , [\gamma , \alpha ]] + 
[ \alpha , [\beta , \gamma ] ] = 0.$$

\noindent \it Proof. \rm\quad To see this at a given point, we need to be slightly more
careful concerning our choice of local holomorphic coordinates, choosing 
\it canonical \rm holomorphic normal coordinates.  
The metric being K\"ahler and real analytic, we can choose local holomorphic
coordinates $z_1, \ldots , z_n$ in a neighbourhood of a given point so that the
K\"ahler form may be written as ${i\over 2} \sum h_{ij}\ dz_i \wedge d \bar z_j$,
where $$ h_{ij} = \delta _{ij}\ + \sum _{k,l}  c_{ijkl} z_k \bar z_l \ + \ O(\, |z|^3 ),$$
with $ c_{ijkl} = {1\over 2} R_{i \bar j k \bar l}$,
where $R$ denotes the Riemannian curvature tensor
(see Appendix 1 from [13], or Exercise 9 on p. 188 of [28]). 
For small $|z|$, the inverse matrix $h^{ij} = 
 \delta _{ij}\ - \sum _{k,l} c_{ijkl} z_k \bar z_l \ + \ O(\, |z|^3 )$.
Now, we use (2.2) to give a formula for $[\alpha , \beta ]$ valid in a neighbourhood of our given point.  Thus we can then use our local 
definition of the bracket to calculate $[\gamma , [\alpha , \beta ]]$ at our 
given point --- although this now involves taking two derivatives, both are 
with respect to the holomorphic rather than anti-holomorphic coordinates.
Since the forms involved are
$$ \theta \ip \beta = {i\over 2} \sum a_{ij} h^{jk} b _{kl} 
 dz_i \wedge d \bar z _l ,$$
and similar terms, when we evaluate at the given point, the derivatives of
$h^{ij}$ do not contribute.  The conclusion therefore is that, by using canonical holomorphic coordinates, we can calculate 
$[\gamma , [\alpha , \beta ]]$ \it at the given point \rm as if the metric were just the flat metric.  The claimed result then reduces to the Jacobi identity on vector fields of type $(1,0)$.

\proclaim Corollary 2.5.  For $\alpha , \beta \in A^{1,1}$,
we have $\partial [\alpha , \beta ] = [\partial \alpha ,
\beta ] - [ \alpha , \partial \beta ]$.\par

\noindent\it Proof.  \rm Applying (2.4),
$$ [\omega , [ \alpha , \beta ]] + [\beta , [\omega , \alpha ]] + 
[ \alpha , [\beta , \omega ] ] = 0,$$ which by (2.1) implies
the claimed result.\bigskip

 Extending these results in a straightforward
 way, we can show that 
$( \bigoplus _{p\ge 0} A^{p,1} , [\ , \ ], 
\partial )$ has the 
structure of a differential graded Lie algebra.

The results of this section are not used in any crucial way in the 
following sections, but they do indicate an underlying structure behind 
the calculations we perform.

\bigskip

\noindent \bf{3.\quad Geodesics}
\medskip\rm

Before proceeding further, we shall need various identities 
on real $(1,1)$-forms 
involving interior and exterior products,  
where at least one of the forms is assumed primitive. 

\proclaim Lemma 3.1.  If $\alpha , \beta$ are real $(1,1)$-forms, with $\alpha = \theta \ip
\omega$ and $\alpha$ primitive, then
$$ (\theta \ip \beta ) \wedge \omega ^{n-1} = - (n-1) \ \alpha \wedge \beta \wedge 
\omega ^{n-2}.$$  In particular, $ (\theta \ip \alpha ) \wedge \omega ^{n-1} 
= - (n-1)\, \alpha ^2 \wedge 
\omega ^{n-2}$.

\noindent\it Proof.\quad\rm  We check the identity at a given point $P$ by
simultaneously diagonalising the forms $\omega$ and $\alpha$, that is choosing
local coordinates $z_1 , \ldots , z_n$  so that 
$\omega = {i\over 2} \sum d z_i \wedge d \bar z_i$ and 
$\alpha = {i\over 2} \sum a_{ii}\, d z_i \wedge d \bar z_i$ at $P$.  
Since $\alpha$ assumed primitive, we have $\sum _i a_{ii}\ =0$ at $P$.  Thus, 
recalling that $({i\over 2})^n dz_1 \wedge d \bar z_1 \wedge \dots \wedge 
dz_n \wedge d \bar z_n = \omega ^n /n!$, at the given point $P$, 
$$ \eqalign{(\theta \ip \beta ) \wedge \omega ^{n-1} 
& = \iha \sum _{i,j} 
a_{ii} b_{ij}\, dz_i \wedge d\bar z_j \wedge \omega ^{n-1}\cr 
& = \big ( \sum _i a_{ii} b_{ii}\big )\ \omega ^n /n \cr
& = - (n-1)\, \alpha \wedge \beta \wedge \omega ^{n-2}}.$$\bigskip

\proclaim Lemma 3.2. (a)\quad Suppose $\alpha , \beta$ are real $(1,1)$-forms, with $\alpha =
\theta \ip \omega$. If $\beta$ is primitive, then
$$2(\theta \ip \beta ) \wedge
\alpha\wedge\omega ^{n-2} = - (n-2)\, \alpha ^2 \wedge \beta \wedge\omega ^{n-3}.$$
If $\alpha$ is primitive, then
$$(\theta \ip \alpha ) \wedge \beta \wedge \omega ^{n-2} + (\theta \ip \beta ) \wedge
\alpha\wedge\omega ^{n-2} = - (n-2)\, \alpha ^2 \wedge \beta \wedge\omega ^{n-3}.$$
In particular, if both $\alpha$ and $\beta$ are primitive, then 
$$(\theta \ip \alpha ) \wedge \beta \wedge \omega ^{n-2} = (\theta \ip \beta ) \wedge
\alpha\wedge\omega ^{n-2} = - \ha (n-2)\, \alpha ^2 \wedge \beta \wedge\omega
^{n-3}.$$  Furthermore, if $\alpha _1, \alpha _2 , \alpha _3$ are all primitive, then $$ (\theta _1 \ip \alpha _2 + \theta _2 \ip \alpha _1 ) \wedge \alpha _3 \wedge \omega ^{n-2}
= - (n-2) \alpha_1 \wedge \alpha _2 \wedge \alpha _3 \wedge \omega ^{n-3}.$$
\noindent (b)\quad   With notation as above, 
suppose $\alpha$ is primitive; then
$$  \beta ^2 \wedge ( \theta \ip \alpha ) \wedge \omega ^{n-3} 
+ 2 \alpha \wedge \beta \wedge ( \theta \ip \beta ) \wedge \omega ^{n-3} 
= -(n-3)\, \alpha ^2 \wedge \beta ^2 \wedge \omega ^{n-4},$$
where the right-hand side is zero for $n\le 3$.  Moreover,
$$ (\theta \ip \theta \ip \beta ) \wedge \beta \wedge 
\omega ^{n-2} + (\theta \ip \beta )\wedge (\theta \ip \beta )
\wedge \omega ^{n-2} = - (n-2)\, \alpha \wedge \beta \wedge 
(\theta \ip \beta ) \wedge \omega ^{n-3}.$$
\medskip

\noindent\it Proof. \rm (a)\quad With coordinates at $P$ chosen as in (3.1), observe
that  $$\alpha^2 \wedge \beta \wedge \omega ^{n-3} = \big ( \sum _{{l,k,i} \dist }
a_{ll} a_{kk} b_{ii} \,
\big )\ \omega ^n \, /n(n-1)(n-2) $$
and 
$$ (\theta \ip \beta ) \wedge \alpha \wedge \omega ^{n-2}  
 = 
\big ( \sum _{k \ne i} a_{ii} b_{ii} a_{kk} \big )\, \omega ^n /n(n-1).$$
If now
$\beta$ is primitive (and so $\sum b_{jj} =0$), then 
$$ \eqalign{ (\theta \ip \beta ) \wedge \alpha \wedge \omega ^{n-2}  
& = \big ( - \sum _{{i,j,k} \dist } a_{ii} b_{jj} a_{kk}
\  - \sum _{k \ne i} a_{ii} b_{kk} a_{kk}\, \big )\ \omega ^n /n(n-1), \cr
& = - (n-2)\, \alpha ^2 \wedge \beta \wedge\omega
^{n-3} - (\theta \ip \beta ) \wedge \alpha \wedge \omega ^{n-2},\cr
}$$ and hence the claim.  If instead $\alpha$ is primitive (and so 
$\sum a_{ll} =0$), then 
$$ \eqalign{ (\theta \ip \beta ) \wedge \alpha \wedge \omega ^{n-2}  
& = \big ( - \sum _{{l,i,k} \dist } a_{ll} b_{ii} a_{kk}
\  - \sum _{k \ne i} a_{kk}^2 b_{ii}\, \big )\ \omega ^n /n(n-1) \cr 
& = - (n-2)\, \alpha ^2 \wedge \beta \wedge\omega
^{n-3} - (\theta \ip \alpha ) \wedge \beta \wedge \omega ^{n-2}.
\cr}$$
If both $\alpha$ and $\beta$ are primitive, the next statement follows immediately.

For the final part, we observe from 
the previous statement  that 
$$ (\theta _1 \pm \theta _2)\ip (\alpha _1 \pm \alpha _2) \wedge \alpha _3 \wedge \omega ^{n-2}
= - \ha (n-2) (\alpha_1 \pm \alpha _2 )^2 \wedge \alpha _3 \wedge \omega ^{n-3}, $$ from which the claim follows.

\noindent \rm (b)\quad With notation as in (a),
$$ \beta ^2 \wedge (\theta \ip \alpha ) \wedge \omega ^{n-3} =
\bigg ( \sum _{{i,j,k} \dist } ( - a_{ii} ^2 b_{jk} b_{kj}  
 + a_{ii}^2 b_{jj} b_{kk} )
\bigg ) \omega ^n /n(n-1)(n-2)
,$$
 the minus sign appearing since $d z_j \wedge d \bar z_k 
\wedge d z_k \wedge d \bar z _j = - d z_j \wedge d \bar z _j 
\wedge d z_k \wedge d \bar z_k$, 
and 
$$ \eqalign{ & 2 \alpha \wedge \beta \wedge (\alpha \ip \beta ) \wedge \omega ^{n-3} = \cr
&  \bigg (  \sum _{{i,j,k} \dist } ( - a_{ii} a_{jj} b_{jk} b_{kj}  
- a_{ii} a_{kk} b_{jk} b_{kj} + a_{ii} a_{jj} b_{jj} b_{kk} + a_{ii} a_{kk} b_{jj} b_{kk}) 
\bigg ) \omega ^n /n(n-1)(n-2).
\cr}$$ 
Using the fact that, for $i,j,k$ distinct, $-a_{ii} - a_{jj} - a_{kk}$ is the sum of the $a_{ll}$ with $l$ distinct from $i,j,k$, we see that the sum of the 
above two forms is therefore 
$$ \bigg (  \sum _{{i,j,k,l} \dist } (  a_{ii} a_{ll} b_{jk} b_{kj} - 
a_{ii} a_{ll} b_{jj} b_{kk}) \bigg ) \omega ^n /n(n-1)(n-2),$$ which then can be identified as $-(n-3) \alpha ^2 \wedge 
\beta ^2 \wedge \omega ^{n-4}.$

The second part of (b) follows from a similar 
(but simpler) local calculation.
\bigskip

\noindent \it Remark 3.3.\rm\quad For $\alpha$ a primitive real
$(1,1)$-form, one can show similarly, using the formula 
for the adjoint Lefshetz operator 
$\Lambda$ on page 114 of [4], that $ \theta \ip \alpha 
= - \ha \Lambda (\alpha \wedge \alpha )$; if both $\alpha$
and $\beta$ are primitive real $(1,1)$-forms, then
$\theta \ip \beta + \phi \ip \alpha = - \Lambda (\alpha \wedge \beta )$.  
When $\alpha$ is also closed (and hence 
harmonic), it follows from the Hodge identity $i \dbar ^* 
= [\Lambda , \d ] $ that 
$$ -i \dbar ^* (\alpha \wedge \alpha ) = -2 \d (\theta \ip \alpha ) = 
[\alpha , \alpha ].$$   So if a primitive 
harmonic real $(1,1)$-form $\alpha$ satisfies 
 $[ \alpha , 
\alpha ] = 0$, or equivalently by (2.3) that 
 $\theta \ip \alpha$ is
closed,  then $ i \dbar ^* (\alpha \wedge \alpha ) =
0$.  The Hodge decomposition corresponding to $\dbar$ then implies that the closed form 
$\alpha \wedge \alpha$ is harmonic, and hence so too is 
$\theta \ip \alpha$.  In the general case, 
given a primitive harmonic real $(1,1)$-form $\alpha$, we 
prove below that $[\alpha ,\alpha ]$ is primitive, and it is 
$\d$-exact by (2.3); thus $[\alpha , \alpha ]$ is
closed (equivalently, $\dbar$-closed) if and only if it is 
zero, which we have just seen 
happens if and only if $\theta \ip \alpha$ is harmonic.
\medskip

\proclaim Lemma 3.4.  If $\alpha$ is  primitive, closed, real
 $(1,1)$-form (hence also harmonic),
then $$\bar \partial ^* (\theta \ip \alpha ) =  i \partial \langle \alpha , \alpha \rangle 
,$$ where $\langle \, ,\, \rangle$ here denotes the pointwise inner-product of forms.  
Equivalently, this says that $[\alpha , \alpha ]$ is
primitive.

\noindent \it Proof. \rm We check this at a given point.  Choose complex 
normal coordinates as in Section 2 (canonical coordinates 
are not needed here).  Given 
$\alpha$, we write $\alpha = {i\over 2} \sum a_{ij} dz_i \wedge d \bar z_j $,
and observe that $$ \langle \alpha , \alpha \rangle  
=  \sum _{i,j} | a_{ij}|^2 + 
O(|z|^2).$$  We now write as before 
$$ \theta \ip \alpha = {i\over 2} \sum _{i,j,k} \tilde a_{ik} a_{kj} \,
dz_i \wedge d \bar z_j = {i\over 2} \sum  _{i,j,k} a_{ik} a_{kj} \,
dz_i \wedge d \bar z_j + O(|z|^2). $$
The formula for $\bar \partial ^* $ on page 113 
of [4]  shows that
$$ \bar \partial ^* (\theta \ip \alpha ) = - * \dbar *  (\theta \ip \alpha ) 
=  i \sum _{i,j,k} {{\partial (a_{ik} a_{kj})} \over {\partial 
z_j }} dz_i + O(|z|).$$
Expanding out, $$ \bar \partial ^* (\theta \ip \alpha ) = 
 i\big ( \sum _{i,j,k} {{\partial a_{ik}} \over {\partial 
z_j }}  a_{kj}\, dz_i + \sum _{i,j,k}  a_{ik} {{\partial a_{kj}} \over {\partial
z_j }}\, dz_i \big ) + O(|z|).$$ The second term in the bracket is 
$ \sum _{i,j,k}  a_{ik} {{\partial a_{jj}} \over {\partial 
z_k }} dz_i $ since $\alpha$ is closed, and this is zero at the point 
since $\sum _j a_{jj} = O(|z|^2)$ from the primitivity of $\alpha$.
The first term can however be rewritten as 
$\sum _{i,j,k} {{\partial a_{jk}} \over {\partial 
z_i }}  a_{kj}\, dz_i $ since $\alpha$ is closed, which in turn can be 
written as $ \partial ( \sum _{j,k} |a_{jk}|^2 )$.  Thus at the point
we have $\bar \partial ^* (\theta \ip \alpha ) = 
 i \partial \langle \alpha , \alpha \rangle $, 
and hence this identity holds everywhere.
This however is equivalent to the condition that $[\alpha ,
\alpha ]$ is primitive, from the fact 
that $\langle \alpha , \alpha \rangle = \Lambda (\theta \ip \alpha)$, 
the Hodge identity $ i \bar \partial ^* = [ \Lambda , \d ]$,
and (2.3).\medskip

Let us now apply (3.4) to give information concerning 
the Hodge decomposition
of $ \theta \ip \alpha$.  We write
$$ (\theta \ip \alpha ) = (\theta \ip \alpha )^h + \dbar \gamma _1 + 
\dbar ^* \gamma _2 .$$  Hence
$$ i \d \langle \alpha , \alpha \rangle = \dbar ^* (\theta \ip \alpha ) = \dbar ^* \dbar 
\gamma _1 $$ for some $(1,0)$-form $\gamma _1$.  This 
latter term  may be written
as $$ \Delta _{\dbar } \,\gamma _1 = \Delta _{\d }\, \gamma _1 = \d \d ^* 
\gamma _1 + \d ^* \d \gamma _1 .$$  
Hodge decomposition for $\d$ now implies that $ \d ^* \d \gamma _1 = 0$, and 
hence $ \d \gamma _1 = 0 $.  Thus $ \dbar \gamma _1 = - i \dbar \d f$ for some
real function $f$.  The Hodge decomposition therefore  reads 
$$ (\theta \ip \alpha ) = (\theta \ip \alpha )^h +  i \d \dbar f + 
\dbar ^* \gamma _2 , \eqno (\dag )$$ 
where  $\dbar ^* \gamma _2 = \d ^* \bar \gamma _2$, 
and so 
$$ (\theta \ip \alpha ) = (\theta \ip \alpha )^h +  i \d \dbar f + i
\d ^* \dbar ^* \Gamma $$ for some $(2,2)$-form $\Gamma$. 

This decomposition has another special property;
 observe that 
 $$ i \d 
\langle \alpha , \alpha \rangle = 
\dbar ^* (\theta \ip \alpha ) = 
i \dbar ^* \d \dbar f =  
- i \d \dbar ^* \dbar f =  - i \d \Delta _{\dbar} f .$$
Thus $\d ( \Delta _{\dbar} f + \langle \alpha , \alpha \rangle ) = 0,$ and so 
$ \Delta _{\dbar} f + \langle \alpha , \alpha \rangle $ is constant on the manifold.  
The constant may be found by integrating over the manifold.  Clearly $
(\Delta _{\dbar} f) \omega ^n$ has integral zero, and the calculation in (3.1) 
shows that $$ \int _X \langle \alpha , \alpha \rangle \, \omega ^n = n 
\int _X (\theta \ip \alpha )\wedge 
\omega ^{n-1} = -n(n-1) \int _X \alpha ^2 \wedge \omega ^{n-2} = n(n-1) A,$$ 
where $A = -  \alpha ^2 \cp \omega ^{n-2} >0$.
Thus $$ -\Delta _{\dbar} f = \langle \alpha , \alpha \rangle - n(n-1)A.$$

We now wedge equation (\dag ) with $\omega ^{n-1}$.  We have seen that 
$(\theta \ip \alpha )\wedge 
\omega ^{n-1} =  \langle \alpha , \alpha \rangle \, \omega ^n /n$; since 
$(\theta \ip \alpha )^h \wedge \omega ^{n-1}$ is harmonic, it is a constant 
multiple of $\omega ^n$, where the constant is clearly $(n-1)A$.  A standard 
calculation shows that $i \d \dbar f \wedge \omega ^{n-1} = 
- (\Delta _{\dbar} f) \, \omega ^n /n $, which is therefore 
$ ({1\over n} \langle \alpha , \alpha \rangle - (n-1)A)\, \omega ^n$.  
We deduce therefore:

\proclaim Proposition 3.5.  In the above decomposition (\dag ), the form
$  (\theta \ip \alpha ) - (\theta \ip \alpha )^h - i \d \dbar f = 
\dbar ^* \gamma _2$ is primitive.

\noindent \it Notation.\rm\quad For any form $\eta $, we shall denote by 
$\eta ^{cl}$ the closed part of the Hodge decomposition of $\eta$, dependent
of course on the metric.  Thus, for $\alpha$ a primitive harmonic real  $(1,1)$-form, equation
(\dag ) says that $ (\theta \ip \alpha )^{cl}
 = (\theta \ip \alpha )^h +  i \d \dbar f ,$ and then (3.5) says that 
$  (\theta \ip \alpha )^{cl} \wedge \omega ^{n-1} = 
(\theta \ip \alpha ) \wedge \omega ^{n-1}.$  Moreover, if we write 
$ (\theta \ip \alpha )^{cl} = (\theta \ip \alpha ) + \gamma $, then 
$\gamma = - i \d ^* \dbar ^* \Gamma$ (for some $(2,2)$-form $\Gamma$) and is 
a primitive element of $d^* A^3$.
\bigskip
We now apply these results to get information about curves 
on the manifold $\Ktilde _1$.
Suppose now that $\omega + \alpha _1 t + 
\ha \alpha _{11} t^2 + 
\si \alpha _{111} t^3 +\ldots$ is a deformation 
of $\omega$ such that 
$$\eqalign{  ( \omega + \alpha _1 t + & \ha \alpha _{11} t^2 + \ldots )^n =
\cr &  \omega ^n + (n \omega ^{n-1}\wedge\alpha _1 ) \, t + \ha \big ( n \omega
^{n-1}\wedge\alpha _{11}  +  n(n-1)\, \omega ^{n-2}\wedge \alpha _1 ^2 \big )\,
t^2 + O(t^3 ).\cr}$$ The reason for the slightly odd subscripts on the coefficients
will become clear in the next section.
Considering the degree one term, we have seen in the Introduction that this implies
that $\alpha _1$ is primitive and closed, and hence harmonic.  If we  set $$\alpha _{11} =  (\theta _1 \ip \alpha _1 )^{cl}  + \xi  ,$$ 
where $\theta _1$ is defined by the equation $\theta _1 \ip \omega = \alpha _1$, then $\xi $ is
a real and closed $(1,1)$-form.  The equation in $t^2$ implies that 
$$  n \, (\theta _1 \ip \alpha _1 )^{cl} \wedge \omega ^{n-1}
+ n\, \xi  \wedge 
\omega ^{n-1} +  n(n-1)\, \omega ^{n-2}\wedge \alpha _1 ^2 = 0, $$
and so 
by (3.1) and (3.5), $\xi$ is primitive, and hence also harmonic.  

\proclaim Theorem 3.6.  Given a 1-parameter family 
on the manifold $\Ktilde _1$ through
$\omega$ in the direction $\alpha_1$,  with $\alpha_1$ primitive harmonic  
 and
$-\int \omega ^{n-2} \wedge \alpha_1^2 = 1$,
we can write it as 
$$\omega (t) = \omega + \alpha_1 t + \ha(  (\theta _1 \ip \alpha _1 )^{cl} 
\,  + \xi  )\, t^2 + O(t^3),$$   where $\xi$ is a primitive 
harmonic form.
  If $\omega (t)$  is
parametrised by arclength, then \break $\int_X \alpha_1 \wedge \xi  \wedge \omega ^{n-2}
= 0$.  If, furthermore, the 1-parameter family is a geodesic, then $\xi  =0$.

\noindent\it Proof.\rm\quad The metric on $\Ktilde _1$ is just the 
Hodge metric.  Since 
$$ {{d \omega}\over{dt}} = \alpha _1 +  \alpha _{11} t + \ha \alpha_{111} t^2 + \ldots
,$$ the condition that this has unit norm is  
$$ - \int _X (\omega + \alpha _1 t + \ha \alpha_{11} t^2 + \ldots )^{n-2} \wedge
(\alpha _1 +  \alpha _{11} t + \ha \alpha_{111} t^2 + \ldots )^2 = 1.$$
Equating terms in $t$ to zero gives
$$ \int_X 2 \alpha_1 \wedge \alpha _{11} \wedge \omega ^{n-2} + (n-2) \int _X \alpha _1
^3 
\wedge \omega ^{n-3} = 0.$$

Recall that $\alpha _{11} = ( \theta_1 \ip \alpha_1 )^{cl} 
\, +  \xi   = 
( \theta_1 \ip \alpha_1 ) - i d^* \dbar ^* \Gamma  
\, +  \xi $.  Observe that $\alpha_1 \wedge \omega ^{n-2}$ is harmonic (and hence the
image of a  harmonic form under the Hodge 
$*$-operator), and so 
$\int _X \alpha_1 \wedge \omega ^{n-2} \wedge d^* \dbar ^* \Gamma 
=0$.
From (3.2)(a), $$2 \alpha _1 \wedge (\theta_1 \ip \alpha_1) \wedge \omega ^{n-2} 
= -(n-2)\, \alpha _1 ^3 
\wedge \omega ^{n-3},$$  and thus the condition obtained from equating terms in $t$ to
zero is  $$ \int_X \alpha_1 \wedge \xi  \wedge \omega ^{n-2} = 0.$$
\medskip

Consider now $\omega (t)$ for $0\le t \le \delta$, with $\delta$ small, and suppose
that it is geodesic; we claim that $\xi  = 0$.  If $ \xi \ne 0$, 
we can choose a family
$\xi (t) = \xi   + t \xi _1 + \dots $ of tangent vectors of the manifold $\Ktilde _1$ 
along the curve $\omega (t)$.  We let $p(t) = t(\delta -t)$,
although the argument also works replacing $t(\delta -t )$ by
other polynomials $p(t)$ vanishing at $0$ and $\delta$,
 and deform the curve
$\omega (t)$ using the  tangent field $p(t) \xi (t)$.  This yields a family of nearby
curves 
$\omega _s (t)$ on $\Ktilde _1$, with the same initial and end points, of the form 
$$\omega _s (t) = \omega (t) + s\, p(t) \,\xi (t) + O(s^2 ).$$
Note that $$\dot \omega _s (t) = \dot\omega (t) + s\, p'(t)\, \xi (t) + 
s \, p(t)\, \dot \xi (t) + O(s^2 ).$$  We use the fact that $\omega (t)$ 
minimizes the energy.

The energy of $\omega _s$ is
$ - \int _0 ^{\delta} \omega _s (t) ^{n-2} \cp \dot\omega _s (t) ^2 \, dt $, which
is then of the form 
$$- \int _0 ^{\delta} \big(\omega (t) + s\, p(t)\, \xi (t) + O(s^2 )\big) ^{n-2}
\cp 
\big( \dot\omega (t) + s\, p'(t)\, \xi (t) + 
s\, p(t)\, \dot \xi (t) + O(s^2 )\big)^2 \, dt,
$$
where as before $\cp$ denotes the cup product of 
cohomology classes represented by closed forms.
If $\omega (t) = \omega _0 (t)$ is geodesic, then
$$ \eqalign{ & (n-2) \int _0 ^{\delta}  p(t) \ \omega (t) ^{n-3}\cp \xi (t)
\cp 
\dot\omega (t) ^2 \, dt \, + \, 2 \int _0 ^{\delta} p'(t) \ \omega (t)
^{n-2}\cp 
\dot\omega (t)\cp \xi (t) \, dt \cr
& + \,  2 \int _0 ^{\delta} p(t) \ \omega (t)
^{n-2}\cp 
\dot\omega (t)\cp \dot\xi (t) \, dt \, =\ 0.
\cr}$$

The term in $\delta ^2$ is $2 \,\omega ^{n-2} \cp 
\alpha _1 \cp \xi  \, \int _0 ^{\delta} p'(t) \,  dt$, which is zero
as required.  Expanding out, the terms in 
$\delta ^3$ are as follows :
$$ \eqalign{ & (n-2)\, \omega ^{n-3} \cp 
\alpha _1 ^2 \cp \xi  \, \int _0 ^{\delta} \{ p(t) + 2t p'(t) \} \,  
dt \cr & + \  2\, \omega ^{n-2} \cp 
\alpha _{11} \cp \xi \, \int _0 ^{\delta} t p'(t) \,  dt \ + \ 
2\, \omega ^{n-2} \cp 
\alpha _1 \cp \xi _1 \, \int _0 ^{\delta} \{ t p'(t) + p(t)\} \, 
dt .\cr }$$

The third term here is clearly zero.  Note that by 
definition of $\alpha _{11}$, 
$$  \omega ^{n-2} \wedge \alpha _{11} \wedge \xi  =  \omega ^{n-2} \wedge 
(\theta _1 \ip \alpha _1 )
\wedge \xi  \,  - \omega ^{n-2} \wedge \xi  \wedge d^* \dbar ^* \Gamma 
+ \,  \omega ^{n-2} \wedge \xi ^2,$$
where by (3.2)(a) 
$$ 2 \, \omega ^{n-2} \wedge (\theta _1 \ip \alpha _1 ) \wedge \xi  \ =\ 
- (n-2) \, \omega ^{n-3} \wedge
\alpha _1 ^2 \wedge \xi  .$$  
Observe also that $\int _X \omega ^{n-2} \wedge \xi  \wedge
d^* \dbar ^* \Gamma = 0$, and so 
$$ 2\, \alpha _{11} \cp \xi \cp  \omega ^{n-2} = - (n-2) \, \omega ^{n-3} \cp
\alpha _1 ^2 \cp \xi    + 2\,  \omega ^{n-2} \cp \xi  ^2.$$
  
Collecting then the terms in $\omega ^{n-3} \cp
\alpha _1 ^2 \cp \xi  $, we obtain a coefficient which is  
$$\int _0 ^{\delta} \{ p(t) + 2t p'(t) - t p'(t) \} \,  
dt \ = \ \int _0 ^{\delta} \{ p(t) + t p'(t) \} \, dt \ = \ 0.$$
  Thus, only one term of order $\delta ^3$ survives, and
this is 
$$ 2 \omega ^{n-2} \cp \xi  ^2 \, \ \int _0 ^{\delta} t p'(t)\, dt \ =\ 
- \th \delta ^3 \, \omega ^{n-2} \cp \xi  ^2 .$$
As this has to be zero,  we deduce 
that $\xi $ is zero (recalling the Hodge Index theorem and that $\xi $ is primitive), and the Theorem is proved.
\bigskip
\noindent \it Remarks 3.7.\quad \rm  
Note that 
there is an equality on harmonic parts $ \alpha _{11} ^h = (\theta _1 \ip
\alpha_1)^h $, and so the expansion for the geodesic $\bar\omega (t)$ on $\calK _1$,
considered as a  submanifold of the harmonic forms $\calH ^{1,1} (\omega )$, is just 
$$ \bar\omega (t) = \omega + \alpha _1 t + \ha (\theta _1 \ip \alpha_1)^h\, t^2 + O(t^3).$$
There is no reason why in general $\alpha _{11}$ should be harmonic.  We note that  
$\alpha _{11} = (\theta _1 \ip \alpha _1)^{cl}$ being  
harmonic is equivalent to
$\dbar ^* (\theta _1 \ip \alpha _1) = 0$,
which by (3.4) is equivalent to  
$\langle \alpha _1 , \alpha _1 \rangle $ being constant over the manifold. 
Thus, if for instance $h^{1,1} > n^2$, at any point of 
the manifold the harmonic real $(1,1)$-forms are linearly dependent, 
and so some non-zero harmonic real $(1,1)$-form $\alpha$ vanishes there.
Since no primitive form is a multiple of $\omega$ at the
point, it follows that $\alpha$ is primitive.  Since  
$\langle \alpha , \alpha \rangle $ cannot be constant on the manifold,
we deduce that 
 $(\theta \ip \alpha )^{cl}$ is not harmonic.  
For example, we can deduce 
 for $h^{1,1} > n^2$ that 
the curved \Kahler cone $\Ktilde$ is not linear, and in 
particular that there exists a primitive harmonic real $(1,1)$-form 
$\alpha$ with $\alpha ^2$ not harmonic; these facts may
be compared with the results from [7].
In the example of  $\omega$ representing the flat metric 
on a complex torus,
the harmonic real $(1,1)$-forms $\alpha$ do  
however have $\langle \alpha ,\alpha \rangle $ constant, with 
the forms $\theta \ip \alpha $
all being harmonic.

\bigskip

\noindent \bf{4.\quad Curvature calculations}
\medskip\rm
In this section, we let $\alpha _1 ,  \ldots , \alpha _r$
denote an orthonormal basis for the tangent space to the manifold $\Ktilde _1$ at
$\omega$, hence a basis of the primitive harmonic real $(1,1)$-forms such that
$$ - \omega ^{n-2} \cp \alpha _i \cp \alpha _j = \delta _{ij}$$ for all $i,j$.
This will determine a local normal coordinate system $(t_1 ,\ldots ,t_r)$; 
locally the
corresponding points of $\Ktilde _1$ will be of the form
$$\omega (t_1 ,\ldots ,t_r) = \omega + \sum _1 ^r \alpha _i t_i + {1\over 2}\sum _{i,j}
\alpha _{ij} t_i t_j + {1\over 6} \sum _{i,j,k}
\alpha _{ijk} t_i t_j t_k + O(t^4).$$ 

The sectional curvature corresponding to the plane spanned by $\alpha _i$ and $\alpha _j$ in
the tangent space of $\Ktilde _1$ at $\omega$ may be determined by expressing the
metric in powers of the normal coordinates $t_i$ and $t_j$ (all the other $t_k$ being held
to be zero).   For simplicity of notation, we shall do this for $\alpha_1$ and $\alpha_2$,
and suppress all the coordinates $t_k$ for $k>2$.  If the metric is written as $\sum
g_{ij} dt_i \otimes dt_j$, then
$g_{12}$ has an expansion
$$ g_{12} (t_1 , t_2 ) = \th R_{1212}\, t_1 t_2 + O(|t|^3) $$
by for instance [19], page 41, where $R_{1212}$ is the sectional curvature we seek.
Note here that the convention we adopt concerning the last two indices of the Riemannian curvature 
tensor differs from that of [19], but coincides with that of [11].  Writing 
$ \omega (t_1 , t_2 ) $ as 
$$ \omega + \alpha _1 t_1 + \alpha _2 t_2 + \ha \alpha _{11} t_1 ^2
+ \alpha _{12} t_1 t_2 + \ha \alpha _{22} t_2 ^2 + \si \alpha _{111} t_1 ^3 + 
\ha \alpha _{112} t_1 ^2 t_2 + \ha \alpha _{122} t_1 t_2 ^2 + \si \alpha _{222} t_2 ^3 
+ \dots ,$$
we have
$$ \partial \omega (t_1 , t_2 ) / \partial t_1 = \alpha _1 + \alpha _{11} t_1 + \alpha _{12} t_2
+ \ha \alpha _{111} t_1 ^2 + \alpha _{112} t_1 t_2 + \ha \alpha _{122} t_2 ^2 + \dots $$
and 
$$ \partial \omega (t_1 , t_2 ) / \partial t_2 = \alpha _2 + \alpha _{12} t_1 + \alpha _{22} t_2
+ \ha \alpha _{112} t_1 ^2 + \alpha _{122} t_1 t_2 + \ha \alpha _{222} t_2 ^2 + \dots .$$
To calculate $g_{12} (t_1 , t_2 )$, we need
$$ - \int_X \omega (t_1 ,t_2 )^{n-2} \wedge \partial \omega (t_1 , t_2 ) / \partial t_1 
\wedge \partial \omega (t_1 , t_2 ) / \partial t_2 .$$  Picking out the terms in $t_1 t_2$, we find
that 
$$\eqalign{ - \th R_{1212} & =  \omega ^{n-2}\cp \alpha_1 \cp \alpha _{122} +
 \omega ^{n-2}\cp \alpha_2 \cp \alpha _{112} + 
 \omega ^{n-2}\cp \alpha_{11} \cp \alpha _{22} \cr & +  
 \omega ^{n-2}\cp \alpha _{12}^2 
+ (n-2) \omega ^{n-3}\cp \alpha_1 ^2 \cp \alpha _{22}
+ (n-2) \omega ^{n-3}\cp \alpha_2 ^2 \cp \alpha _{11} \cr & +
3(n-2) \omega ^{n-3}\cp \alpha_1 \cp \alpha_2 \cp \alpha _{12} +
(n-2)(n-3) \omega ^{n-4}\cp \alpha_1 ^2 \cp \alpha_2 ^2 .\cr
}$$

We can also calculate the curvature from the fact that 
$g_{11}$ has an expansion
$$ g_{11} (t_1 , t_2 ) = 1 + \th R_{1221}\, t_2 ^2 + O(|t|^3) .$$
To calculate $g_{11} (t_1 , t_2 )$, we need
$$ - \int_X \omega (t_1 ,t_2 )^{n-2} \wedge \partial \omega (t_1 , t_2 ) / \partial t_1 
\wedge \partial \omega (t_1 , t_2 ) / \partial t_1 .$$  
Picking out the terms in $t_2 ^2$,
we find that
$$\eqalign{ \th R_{1212} & =  \omega ^{n-2}\cp \alpha_1 \cp \alpha _{122} 
  +   \omega ^{n-2}\cp \alpha _{12}^2 \cr &
+ \ha (n-2) \omega ^{n-3}\cp \alpha_1 ^2 \cp \alpha _{22}
 + 2(n-2) \omega ^{n-3}\cp \alpha_1 \cp \alpha_2 \cp \alpha _{12} \cr & +
\ha (n-2)(n-3) \omega ^{n-4}\cp \alpha_1 ^2 \cp \alpha_2 ^2 .\cr
}$$
 By symmetry therefore,
$$\eqalign{ \tth R_{1212} & =  \omega ^{n-2}\cp \alpha_1 \cp \alpha _{122} +
 \omega ^{n-2}\cp \alpha_2 \cp \alpha _{112}  \cr & +  
 2 \omega ^{n-2}\cp \alpha _{12}^2 
+ \ha (n-2) \omega ^{n-3}\cp \alpha_1 ^2 \cp \alpha _{22}
+ \ha (n-2) \omega ^{n-3}\cp \alpha_2 ^2 \cp \alpha _{11} \cr & +
4(n-2) \omega ^{n-3}\cp \alpha_1 \cp \alpha_2 \cp \alpha _{12} +
(n-2)(n-3) \omega ^{n-4}\cp \alpha_1 ^2 \cp \alpha_2 ^2 .\cr
}$$
Subtracting our previous expression for $-\th R_{1212}$ from this, we obtain
$$\eqalign{ R_{1212} & =    
\omega ^{n-2}\cp \alpha _{12}^2  - \omega ^{n-2}\cp \alpha _{11} \cp \alpha _{22} \cr &
- \ha (n-2) \omega ^{n-3}\cp \alpha_1 ^2 \cp \alpha _{22}
- \ha (n-2) \omega ^{n-3}\cp \alpha_2 ^2 \cp \alpha _{11} \cr & +
(n-2) \omega ^{n-3}\cp \alpha_1 \cp \alpha_2 \cp \alpha _{12}  .\cr
}$$

Thus far, we have yet to use the information gleaned from the previous section.
 Extending our previous notation (Theorem 3.6) in an obvious way, we know that
$\alpha _{ii} = 
(\theta _i \ip \alpha _i )^{cl} = \theta _i \ip \alpha _i + \gamma_{ii}$, where 
$\gamma _{ii} \in d^* A^3 $ and is primitive, and 
as usual $\theta _i \ip \omega = \alpha _i$.  Moreover, $\gamma _{ii}$ may be recovered from $\alpha _i \wedge \alpha _i$ as follows:  There exists a 
(unique) $\d ^*$-coexact $(1,2)$-form $\eta _i$ such that 
$\alpha _i \wedge \alpha _i = (\alpha _i \wedge \alpha _i )^h + \d \eta _i ;$
then $\theta _i \ip \alpha _i = -\ha \Lambda (\alpha _i \wedge \alpha _i )^h 
- \ha \Lambda \d  \eta _i $.  Using the Hodge identity $[\Lambda , \d ] = i 
\dbar ^*$, we deduce that $\gamma _{ii} = - {i\over 2} \dbar ^* \eta _i \in \d ^* \dbar ^* A^{2,2}$.

We can however also identify 
the cross-terms $\alpha _{12}$, as we know that the curve given by
$t_1 = t_2 = t/\sqrt{2}$ and $t_k =0$ otherwise, will also be a geodesic, corresponding
to the unit tangent  vector $ (\alpha _1 + \alpha _2 )/\sqrt{2} $.  
From this it follows that 
$$ \ha  \alpha _{11} + \ha \alpha _{22} +
\alpha _{12} = \ha ((\theta _1 + \theta _2) \ip 
(\alpha _1 + \alpha _2 ))^{cl}  .$$    Thus 
$$ \alpha _{12} = \ha \big ( \theta _1 \ip \alpha _2 + \theta _2 \ip \alpha
_1 \big ) + \gamma _{12}
=  \Re\, \big ( \theta _1 \ip \alpha _2  \big ) + \gamma _{12} 
=  \Re\, \big ( \theta _1 \ip \alpha _2  \big )^{cl},$$
with  $ \gamma _{12} \in d^* A^3 $ and primitive.
In particular, 
 $$\omega ^{n-1} \cp \alpha _{12} = \int _X \omega ^{n-1} \wedge \Re \, (\theta _1 \ip 
\alpha _2) = - (n-1) \, \alpha _1 \cp \alpha _2 \cp \omega ^{n-2},$$ using (3.1).
By assumption $ -\omega ^{n-2}\cp \alpha _1 \cp \alpha _2 =0$, and thus 
$\omega ^{n-1}\cp \alpha _{12} = 0$.

\proclaim Theorem 4.1. Given primitive harmonic forms $\alpha _1 , \alpha _2$ as above 
on the compact \Kahler manifold $X$, the corresponding sectional curvature of $\calK _1$ is 
given by the formula 
$$ R_{1212} = \int _X (\Im (\theta _1 \ip \alpha _2 ) )^2 \wedge \omega ^{n-2} 
 + \int _X \gamma _{12}^2 \wedge \omega ^{n-2} -
\int _X \gamma_{11} \wedge \gamma _{22} 
\wedge \omega ^{n-2} ,$$
where the first two terms are non-positive.\par

\noindent\it Proof. \quad\rm  From (3.2)(a), we know that for any  
primitive form $\alpha _3$, 
$$ (\theta _1 \ip \alpha _2 + \theta _2 \ip \alpha _1 ) \wedge \alpha _3 \wedge \omega ^{n-2}
= - (n-2) \alpha_1 \wedge \alpha _2 \wedge \alpha _3 \wedge \omega ^{n-3}.$$
We have  
$$\alpha _{12} = \ha (\theta _1 \ip \alpha _2 + \theta_2 \ip \alpha_1 ) 
+ \gamma _{12},$$ with $\gamma _{12}$ 
a primitive real $(1,1)$-form.  Thus,  
$$\eqalign{  \alpha _{12} ^2 \wedge \omega ^{n-2} & = \qu (\theta _1 \ip \alpha _2 + \theta_2 \ip \alpha_1 )^2 \wedge \omega ^{n-2} + 
(\theta _1 \ip \alpha _2 + \theta_2 \ip \alpha_1 )\wedge \gamma _{12} 
 \wedge \omega ^{n-2} + \gamma _{12}^2 \wedge \omega ^{n-2} \cr
& = \qu (\theta _1 \ip \alpha _2 + \theta_2 \ip \alpha_1 )^2 \wedge \omega ^{n-2}  - (n-2) \alpha _1 \wedge \alpha _2 \wedge \gamma _{12} 
 \wedge \omega ^{n-2} + \gamma _{12}^2 \wedge \omega ^{n-2}.\cr}$$

Therefore, the sum 
$$\eqalign{ & \alpha _{12} ^2 \wedge \omega ^{n-2} + 
(n-2) \alpha_1 \wedge \alpha_2 \wedge \alpha _{12} \wedge  \omega ^{n-3} = \cr & 
\qu (\theta _1 \ip \alpha _2 + \theta_2 \ip \alpha_1 )^2 \wedge \omega ^{n-2} + \gamma _{12}^2 \wedge \omega ^{n-2}  + 
(n-2) \alpha_1 \wedge \alpha_2 \wedge (\alpha _{12} - \gamma_{12}) \wedge  \omega ^{n-3} = \cr
& \qu (\theta _1 \ip \alpha _2 + \theta_2 \ip \alpha_1 )^2 
\wedge \omega ^{n-2} + \gamma _{12}^2 \wedge \omega ^{n-2} + \ha (n-2) \alpha _1 \wedge \alpha _2 
\wedge (\theta _1 \ip \alpha _2 + \theta_2 \ip \alpha_1 ) \wedge \omega ^{n-3}
.\cr }$$

This may be written as 
$$ \eqalign{ & - \qu  (\theta _1 \ip \alpha _2 - \theta _2 \ip \alpha _1 )^2 \wedge 
\omega ^{n-2} + \gamma _{12}^2 \wedge \omega ^{n-2}
+ \ha (\theta _1 \ip \alpha _2 )^2 \wedge \omega ^{n-2} 
 + \ha (\theta _2 \ip \alpha _1 )^2 \wedge \omega ^{n-2} \cr 
& + \ha (n-2)  \alpha _1 \wedge \alpha _2 \wedge 
(\theta _1 \ip \alpha _2 ) \wedge \omega ^{n-3}
+ \ha (n-2)  \alpha _1 \wedge \alpha _2 \wedge 
(\theta _2 \ip \alpha _1 ) \wedge \omega ^{n-3}.\cr }$$
Applying the second identity of (3.2)(b) twice, this 
may be rewritten as 
$$  (\Im (\theta _1 \ip \alpha _2 ))^2 \wedge \omega ^{n-2} 
 + \gamma _{12}^2 \wedge \omega ^{n-2} 
- \ha  (\theta _1 \ip \theta _1 \ip \alpha _2 )\wedge \alpha _2 \wedge \omega ^{n-2}
- \ha  (\theta _2 \ip \theta _2 \ip \alpha _1 )\wedge \alpha _1 \wedge \omega ^{n-2}.$$

A straightforward local calculation verifies that
$$ (\theta _1 \ip \theta _1 \ip \alpha _2 )\wedge \alpha _2 \wedge \omega ^{n-2}
= (\theta _2 \ip \theta _1 \ip \alpha _1 )\wedge \alpha _2 \wedge \omega ^{n-2} ,$$
and hence also the corresponding formula with indices exchanged. 
Hence $$\eqalign{  \alpha _{12} ^2 \wedge \omega ^{n-2} + &
(n-2) \alpha_1 \wedge \alpha_2 \wedge \alpha _{12} \wedge  \omega ^{n-3}  
 = 
(\Im (\theta _1 \ip \alpha _2 ))^2 \wedge \omega ^{n-2} 
 + \gamma _{12}^2 \wedge \omega ^{n-2} \cr &
- \ha  (\theta _2 \ip \theta _1 \ip \alpha _1 )\wedge \alpha _2 \wedge \omega ^{n-2}
- \ha  (\theta _1 \ip \theta _2 \ip \alpha _2 )\wedge \alpha _1 \wedge \omega ^{n-2},\cr}$$ which using (3.2)(a)  
 may be written as 
$$\eqalign{ & (\Im (\theta _1 \ip \alpha _2 ))^2 \wedge \omega ^{n-2} 
 + \gamma _{12}^2 \wedge \omega ^{n-2} + (\theta_1 \ip \alpha _1 )\wedge 
(\theta _2 \ip \alpha _2 )\wedge \omega ^{n-2} \cr &
+ \ha (n-2) \alpha _1 ^2 \wedge (\theta_2 \ip \alpha _2 )\wedge \omega ^{n-3} 
+ \ha (n-2) \alpha _2 ^2 \wedge (\theta_1 \ip \alpha _1 )\wedge \omega ^{n-3}.
\cr}$$
Substituting this into our formula for $R_{1212}$, we get that $R_{1212}$ is the  integral of 
$$\eqalign{  & (\Im (\theta _1 \ip \alpha _2 ))^2 \wedge \omega ^{n-2} 
 + \gamma _{12}^2 \wedge \omega ^{n-2} - \gamma_{11}\wedge (\theta _2 \ip 
\alpha _2 )\wedge \omega ^{n-2} - \gamma_{22}\wedge (\theta _1 \ip 
\alpha _1 )\wedge \omega ^{n-2} \cr & -\gamma_{11} \wedge \gamma _{22} 
\wedge \omega ^{n-2} - \ha (n-2) \alpha _1 ^2 \wedge \gamma _{22} 
\wedge \omega ^{n-2}  - \ha (n-2) \alpha _2 ^2 \wedge \gamma _{11} 
\wedge \omega ^{n-2} ,\cr}$$ which, using 
the third part of (3.2)(a), and $\gamma_{11}, \gamma _{22}$ being primitive, simplifies to 
$$ (\Im (\theta _1 \ip \alpha _2 ) )^2 \wedge \omega ^{n-2} 
 + \gamma _{12}^2 \wedge \omega ^{n-2} -\gamma_{11} \wedge \gamma _{22} 
\wedge \omega ^{n-2} .$$

Finally, we observe that $\Im (\theta _1 \ip \alpha _2 ) = 
- {i\over 2} ( \theta _1 \ip \alpha _2 - 
\theta _2 \ip \alpha _1 )$ is primitive by (3.1), and 
$\gamma _{12}$ is primitive by (3.5), and so the first two 
terms in our formula are non-positive.\medskip

\noindent \it Remark.\quad\rm Our formula (4.1) for 
the sectional curvature may be rewritten as 
$$ n(n-1) R_{1212} = - (\Im (\theta _1 \ip \alpha _2 ), 
\Im (\theta _1 \ip \alpha _2 ) ) - (\gamma _{12} , \gamma _{12}) + (\gamma _{11}, \gamma _{22} ),$$
where the bracket $(\, , \, )$ 
here denotes the global inner-product.  
Compare this formula with the Gauss formula 
for the curvature of a submanifold,
page 48 of [28], relevant because we can consider 
$\tilde \calK _1$ as a submanifold of the infinite dimensional Riemannian manifold of all real 
$(1,1)$-forms with top power being $\omega ^n$.  The $\gamma _{ij}$ are 
 analogous to the values of an appropriate multiple of
the second fundamental form.
\bigskip
\noindent\it Example 4.2.\quad\rm  An instructive example 
 is provided 
when
$X$ is an $n$-dimensional complex torus.
After choosing a basis for the global holomorphic 1-forms, 
$H^{1,1}(X, \R )$ may be identified as
the hermitian $n \times n$ matrices.
The \Kahler cone corresponds to the positive definite
hermitian matrices, and the form of degree $n$ coming from cup product is some positive multiple
of the determinant.  Thus the normalised \Kahler moduli space
$\calK _1$ may be identified as the space of positive definite
hermitian matrices of determinant 1, that is, the symmetric space $ SL(n, \C )/SU(n)$.  
The Hodge metric on $\calK _1$ is invariant under the action of $ SL(n, \C )$,
corresponding to different choices of bases for the holomorphic 1-forms (not
changing the volume form).  Since
$ SL(n, \C )/SU(n)$ is an irreducible symmetric space, the metric we have defined 
must be a constant multiple of the symmetric space metric; thus in this case $\calK _1$ is  
also complete.  A very concrete description of the symmetric space $SL(n,\R )/SO(n)$ is given in Section
5.4 of [8]; the theory for $ SL(n, \C )/SU(n)$ is entirely analogous, just replacing
 the transpose of a matrix by its 
hermitian conjugate, 
and with the properties symmetric, respectively antisymmetric,  being replaced by hermitian,
respectively skew-hermitian.
 The tangent space at $I_n$ is given 
by the trace-free hermitian matrices; given two such matrices $A$ and $B$ representing orthonormal
tangent vectors, 
the sectional curvature of the corresponding tangent plane is ${{\cal B}} ([A,B] , [A,B])$
(cf. [8], Corollary 5.4.1), where the bracket
operation $[ \, , \, ]$ 
is just the commutator of the matrices and ${{\cal B}}$ denotes the standard Killing 
bilinear form.  This then corresponds to the formula we derived above,
modulo an expected constant factor, once we have lifted the \Kahler class to the form 
representing the flat metric --- here of course $\gamma _{11}, \gamma _{22}$ and $\gamma_{12}$ are all zero.

For a trace-free skew-hermitian matrix $C = (c_{ij})$, we have that
$${{\cal B}}(C, C) = 2n \ \tr (C^2) = - 2n\  \tr (CC^*) = -2n \sum _{i,j} |c_{ij}|^2.$$
By choosing a suitable basis for the 1-forms, we may assume that the \Kahler form corresponds
to the identity matrix and that 
$A$ is diagonal, and then a routine 
check verifies that 
${{\cal B}}([A,B] , [A,B]) = -2n \sum _{i,j} (a_{ii} - a_{jj})^2 |b_{ij}|^2$.   We
observe in passing that, when $n\ge 3$, the curvature in this example is not \it strictly \rm negative.
\bigskip
There is however a second expression for the sectional curvature, which is
relevant when considering the possible existence of a \it lower \rm bound
for the sectional curvatures.

\proclaim Theorem 4.3. With the same notation as in (4.1), the 
 sectional curvature $ R_{1212}$ is also 
given by the formula 
$$-3 \alpha _{12} ^2 \cp \omega ^{n-2} + \ha (n-2)(n-3)
\alpha _1 ^2 \cp \alpha _2 ^2 \cp \omega ^{n-4}
 -2 \int _X \gamma _{12}^2 \wedge \omega ^{n-2} 
-\int _X \gamma_{11} \wedge \gamma _{22} \wedge \omega ^{n-2} 
- \ha n(n-1),$$
where here the first three terms are non-negative.\par

\noindent\it Proof. \quad\rm  
We observe first that 
$$ \alpha _{12} ^2 \cp \omega ^{n-2} = - \ha (n-2) \, \alpha _1 \cp \alpha _2 \cp \alpha _{12}
\cp \omega ^{n-3}.$$
To prove this, note that $$\alpha _{12} ^2 \cp \omega ^{n-2} = \alpha _{12} \cp \alpha _{12} ^h 
\cp \omega ^{n-2} = \ha \int _X (\theta _1 \ip \alpha _2 + \theta _2 \ip \alpha _1 )
\wedge \alpha _{12} ^h \wedge \omega ^{n-2},
$$ and that $\alpha _{12} ^h $ is also primitive, 
and so the claim follows from (3.2)(a).

A further identity is that 
$$ \alpha _{11} \cp \alpha _{22} \cp \omega ^{n-2} =
 - \ha (n-2) \alpha _2 ^2 \cp \alpha _{11} \cp \omega ^{n-3} + \ha n(n-1) .$$
This is 
proved by considering the harmonic form $\alpha _{11} ^h 
-(n-1)\omega$, which is primitive by (3.1), and applying (3.2)(a).
Thus for instance one can in general write  
$$\eqalign{R_{1212} &= - \alpha _{12} ^2 \cp \omega ^{n-2} + \alpha _{11} \cp \alpha _{22} \cp
\omega ^{n-2} - n(n-1) \cr & = 
- \alpha _{12} ^2 \cp \omega ^{n-2} -  \ha (n-2) \alpha _1 ^2 \cp \alpha _{22} \cp \omega
^{n-3} - \ha n(n-1).\cr}$$  

From the first part of (3.2)(b), we have 
$$  \alpha _1 ^2 \wedge ( \theta _2 \ip \alpha _2 ) \wedge \omega ^{n-3} 
+ 2 \alpha _1 \wedge \alpha _2 \wedge ( \theta _2 \ip \alpha _1 ) 
\wedge \omega ^{n-3} 
= -(n-3)\, \alpha _1 ^2 \wedge \alpha _2 ^2 \wedge \omega ^{n-4}.$$
An easy local calculation verifies that 
$$ \alpha _1 \wedge \alpha _2 \wedge ( \theta _2 \ip \alpha _1 ) 
\wedge \omega ^{n-3} = 
\alpha _1 \wedge \alpha _2 \wedge ( \theta _1 \ip \alpha _2 ) \wedge \omega ^{n-3} ,$$
and so 
$$  \alpha _1 ^2 \wedge ( \theta _2 \ip \alpha _2 ) \wedge \omega ^{n-3} 
+ \alpha _1 \wedge \alpha _2 \wedge ( \theta _1 \ip \alpha _2 + 
\theta _2 \ip \alpha _1 ) 
\wedge \omega ^{n-3} 
= -(n-3)\, \alpha _1 ^2 \wedge \alpha _2 ^2 \wedge \omega ^{n-4}.$$
Thus 
$$ \eqalign{  - \ha (n-2) \, \alpha _1 ^2 \cp \alpha _{22} \cp & \omega ^{n-3}  
 +  \ha (n-2) \int _X \alpha _1 ^2 \wedge \gamma _{22} \wedge \omega ^{n-3} 
- (n-2) \, \alpha _1 \cp \alpha _2 \cp \alpha _{12} 
\cp \omega ^{n-3} \cr & + (n-2) \int _X
\alpha _1 \wedge \alpha _2 \wedge \gamma _{12} \wedge \omega ^{n-3} 
= \ha (n-2)(n-3) \, \alpha _1 ^2 \cp \alpha _2 ^2 \cp \omega ^{n-4}.\cr }$$

Using the fact that for any closed real $(1,1)$-form $\alpha$ 
and primitive real $(1,1)$-form 
$\gamma \in d^* A^3$, we have $\int _X \alpha \wedge \gamma \wedge 
\omega ^{n-2} = 0$,
we deduce, using the third and fourth equations from 
(3.2)(a), and 
also the first identity of the current proof, that 
$$ \eqalign{ 
\ha (n & -2) (n-3) \, \alpha _1 ^2 \cp  \alpha _2 ^2 \cp \omega ^{n-4} = \cr 
&  - \ha (n-2) \, \alpha _1 ^2 \cp \alpha _{22} \cp  \omega ^{n-3}  
 +  \int _X \gamma _{11} \wedge \gamma _{22} \wedge \omega ^{n-2} 
\,  + \, 2 \,  \alpha _{12} ^2 
\cp \omega ^{n-2}  +  2 \int _X
\gamma _{12} ^2 \wedge \omega ^{n-2} 
 .\cr }$$
Using this to substitute for 
$ - \ha (n-2) \, \alpha _1 ^2 \cp \alpha _{22} \cp  \omega ^{n-3} $
 into the second expression given above 
for $R_{1212}$, we obtain the formula claimed.

The first and the third terms of this formula are non-negative from the primitivity of the
forms $\alpha ^h _{12}$ and $\gamma _{12}$, and the second term is non-negative 
from the primitivity of the 4-form $(\alpha _1 \wedge \alpha _2 )^h$ and the \
Hodge--Riemann bilinear relations.\bigskip

\noindent \it Discussion 4.4. \rm\quad  We have seen that it is natural therefore extend the question concerning 
semi-negativity of \Kahler moduli curvature to one 
which asks for which \Kahler manifolds the 
sectional curvatures of $\calK _1$ lie between $-\ha n(n-1)$ and $0$.  The 
message of (4.1) and (4.2) is that for both bounds, the term 
$ - \int _X \gamma _{11} \wedge \gamma _{22} \wedge \omega ^{n-2}$ is crucial.
If for instance 
 $\alpha _1 ^2$ is harmonic, then this term is zero and the 
sectional curvature $R_{1212}$ 
has the conjectured bounds.  We deduce immediately 
that if the squares  of harmonic real $(1,1)$-forms on 
$X$ are harmonic, for instance if $X$ is a hermitian symmetric space of compact type, then the bounds are as stated for all the sectional 
curvatures.  Moreover, in Example 4.2, both bounds are in fact achieved when $n>2$.  The formulae (4.1) and (4.3) 
 are however rather stronger than this. 

In the \CY case, we argued in Section 1 that we can reduce down to the case
of simply connected \CY manifolds with $h^{2,0} = 0$, and that in this case 
$\calK _1$ was essentially independent of the complex structure. 
Using an argument via Yau's construction of the Ricci flat metric and Moser's Theorem, we may even take the symplectic form $\omega$ to be fixed as we 
vary the complex structure.  In the general case, we can also fix $\omega$, 
and then vary the complex structure among those compatible with $\omega$ and 
which induce the same Hodge structure on $H^2 (X, \R )$, this latter condition
holding automatically when $h^{2,0} =0$.  In this 
more general case also,  
on the complement of countably many subvarieties in the complex structure moduli space, the normalised \Kahler 
moduli space 
$\calK _1$ is essentially independent of the complex structure.  This 
then leads us to the following criterion:

\noindent \it Criterion.\rm\quad For a given \Kahler class $\bar \omega$ and 
orthonormal primitive classes $\bar \alpha _1$ and $\bar \alpha _2$ in 
$H^{1,1} (X, \R )$, we lift $\bar \omega$ to a \Kahler form $\omega$, which we 
now fix.  We vary the complex structure as described above, and for 
each such compatible complex structure, we lift the classes 
$\bar \alpha _1$ and $\bar \alpha _2$ to their harmonic representatives.
This then gives rise to the forms $\gamma _{11}$ and 
$\gamma _{22}$, and hence the quantity 
 $ \int _X \gamma _{11}\wedge \gamma _{22} \wedge 
\omega ^{n-2}$, 
which we may regard as a function of the complex 
structure.  It follows from $\int _X (\theta _1 \ip 
\alpha _1 )\wedge (\theta _2 \ip \alpha _2) \wedge 
\omega ^{n-2} \ge 0$ that 
 $ \int _X \gamma _{11}\wedge \gamma _{22} \wedge 
\omega ^{n-2}$ is 
bounded below by 
$ - \alpha _{11}\cp \alpha _{22} \cp \omega ^{n-2}$, 
independent
of the complex structure.  If now for some limit point of complex moduli, 
either in the interior or more likely maybe on the boundary, 
the above quantity tends to 0 
as we approach the limit point suitably, then the corresponding sectional 
curvature $R_{1212}$ of $\calK _1$ satisfies 
$-\ha n(n-1) \le R_{1212} \le 0$.  

If we write $\alpha _i \wedge \alpha _i = (\alpha _i \wedge \alpha _i )^h + 
\partial \eta _i$, with $\eta _i$ assumed to be 
$\d ^*$-coexact, then as noted above $\gamma _{ii} = - \ha i \dbar ^* \eta _i$.
Since $\dbar\d \eta _i =0$, we have $\d\dbar \eta _i =0$.  But 
$\d ^* (\dbar \eta _i ) = - \dbar \d ^* \eta _i = 0$, and so 
$\dbar \eta _i$ is harmonic; in particular $\dbar ^* \dbar 
\eta _i = 0$.
Now since $\gamma _{11}$ and $\gamma _{22}$ are primitive, 
$$n(n-1) \int _X \gamma _{11} \wedge \gamma _{22} \wedge 
\omega ^{n-2} = - ( \gamma _{11} , \gamma _{22} ) 
= - \qu ( \dbar ^* \eta _1 , \dbar ^* \eta _2 ) 
=  - \qu (\dbar \dbar ^* \eta _1 ,  \eta _2 ).$$
Since  $\dbar ^* \dbar \eta _1 = 0$, this may be written as 
$$  - \qu (\Delta _{\dbar} \eta _1 ,  \eta _2 ) = 
- \qu (\Delta _\d \eta _1 ,  \eta _2 ) =
 - \qu (\d ^* \d \eta _1 ,  \eta _2 ) 
=  - \qu (\d \eta _1 , \d \eta _2 ).$$
Thus $$ (\alpha _1 \wedge \alpha _1 , \alpha _2 \wedge \alpha _2 )
= ((\alpha _1 \wedge \alpha _1)^h  ,   (\alpha _2 \wedge \alpha _2)^h ) \ + 
\ (\d \eta _1 , \d \eta _2),$$ the first term of which is invariant 
and the second term being (up to the constant) the quantity we wish to know 
about.

When dealing with \CY manifolds with $h^{2,0} = 0$, it is tempting to believe
that one should consider a large complex structure limit point for the above purposes.  The author hopes to return to 
this question in a later paper.

\bigskip

\noindent \bf{5.\quad \Kahler threefolds with $b_2 = 3$} \medskip \rm
\medskip \rm  
Given a form of degree $n$ on $\R ^r$, we define the index cone $W$, as in the 
Introduction, to be the set of points at which the form is positive and the signature
of the associated quadratic form is $(1,r-1)$.
In principle,  one can then
calculate the sectional curvatures of the level set $W_1$.  
The easiest case will be that of cubic forms on $\R ^3$, of relevance to
\Kahler threefolds with $h^{1,1} = 3$. Here we shall have  
only one sectional curvature to consider.

We denote the cubic form by $F$, 
and its Hessian by $H$, also a cubic form.  Choosing a basis for the vector space, these
may be considered as homogeneous cubics in three variables.  We let
$S$ denote the degree 4 invariant of $F$, one of the two basic invariants $S$ and $T$ found
by Aronhold, and written down explicitly in [20]. 
If the cubic is non-singular
and in Weierstrass normal form, these
are multiples of the better known invariants $g_2, g_3$ of the corresponding elliptic
curve,  with in particular  
$ S = 4g_2 /27$.  For a general ternary cubic, $S$ has 25 terms
and $T$ has 103 terms.
 Modulo a positive constant, the discriminant 
$\Delta$ of the cubic is $64 S^3 - T^2$.  Recall that the cubic is singular if and only if
$\Delta = 0$, and in the smooth case the real curve has one or two components, dependent
upon whether $\Delta <0$ or $\Delta >0$.  Thus if $S <0$, we can only have one component,
whilst if $S>0$ we may have one or two components.  
We let $H(X_0 , X_1 , X_2)$ denote the Hessian of
$F$, also a cubic form.

\proclaim Theorem 5.1.  If $F$ is a ternary cubic form with Hessian $H$ and basic degree 4 
invariant $S$, then 
for any point of its index cone $W$, 
the curvature at the corresponding point of the level set $W_1$
is given by the formula
$$ R = - {9\over 4} + {{6^6 S \, F^2} \over {4 H^2 }}.$$\par

\noindent\it Remarks.\quad\rm If in fact we are at a point of $W_1$, then  
by definition $F$ takes value one.  I prefer however to write the formula in the way given, 
since then all the required invariance under changes of coordinates and scalings hold good
and are clear.  Motivated by this formula, Totaro has recently demonstrated a natural extension of it, 
in a similar shape, for forms $F$ of arbitrary degree $d>2$.  As the metric on $W_1$ is given by the restriction of $g_{ij} = - \si \partial ^2 F / \partial X_i \partial X_j$, the above formula may be rewritten as $ R = - {9\over 4} + 
\qu S F^2/h^2 $, where $h = \det (g_{ij})$.

\noindent\it Proof.\quad\rm  Suppose the relevant point of $W_1$ is represented by an
element $L_0$ of $\R ^3$, at which the cubic form  takes the value one.  We can then find 
a real basis 
$\{ L_0 , L_1 , L_2 \}$, with $L_0 ^2 \,\cp L_i = 0$  for
$i=1,2$,  where $\cp$ is given by polarising the cubic form.
We shall let $L_0 \, \cp L_1 ^2 = - p <0$, $L_0 \,\cp L_1 \, \cp L_2 = - q$, 
$L_0 \,\cp L_2 ^2 = - r <0$, and $L_1 ^i \,\cp L_2 ^j = a_{ij}$ for all non-negative 
$i+j = 3$.  One could of course also choose $L_1 ,L_2$ so that $p=r= 1$ and $q=0$,
but that does not significantly simplify the problem and tends to obscure some of
the structure.  In these coordinates, the value of the cubic form at 
$ x_0 L_0 + x_1 L_1 + x_2 L_2$ is
$$ (x_0 L_0 + x_1 L_1 + x_2 L_2)^3 =  x_0^3 - 3 x_0 (p x_1 ^2 + 2q x_1 x_2 + r x_2 ^2)
+ a_{30} x_1 ^3 + 3 a_{21} x_1 ^2 x_2 + 3 a_{12} x_1  x_2 ^2 + a_{03} x_2 ^3 ,$$
and so the level set $W_1$ has equation
$$  F(x_0 ,x_1, x_2) = x_0^3 - 3 x_0 (p x_1 ^2 + 2q x_1 x_2 + r x_2 ^2)
+ a_{30} x_1 ^3 + 3 a_{21} x_1 ^2 x_2 + 3 a_{12} x_1  x_2 ^2 + a_{03} x_2 ^3 = 1.$$

The Hodge metric at a point $(x_0 , x_1 , x_2 )$ then takes the form
$$\eqalign{& - (x_0 L_0 + x_1 L_1 + x_2 L_2)\,\cp (L_0 dx_0 +  L_1 dx_1 +  L_2 dx_2 )^2 \cr
& = - dx_0 ^2 + x_0 (p dx_1 ^2 + 2q dx_1 dx_2 + r dx_2 ^2 ) 
+ 2  (px_1 dx_1 + q (x_1 dx_2 + x_2 dx_1) + r x_2 dx_2 ) dx_0 \cr
& - (a_{30} x_1 + a_{21} x_2) d x_1 ^2 - 2 (a_{21}x_1 + a_{12} x_2 ) dx_1 dx_2 
-(a_{12} x_1 + a _{03} x_2 ) dx_2 ^2 .
\cr }$$

Now set $x_0 = 1 + y_0$, so that
$$ y_0 = (p x_1 ^2 + 2q x_1 x_2 + r x_2 ^2) + O(|x|^3),$$
where $|x|^2 = x_1 ^2 + x_2 ^2$.  Therefore
$$ dy_0 = d x_0 = 2 (p x_1 d x_1 + q x_1 d x_2 + q x_2 d x_1 + r x_2 dx_2) + O(|x|^2).$$

Substituting in, we find that the metric takes the form locally 
$ g = \sum g_{ij} dx _i d x_j ,$ where
$$ g_{11} = - a_{30} x_1 - a_{21} x_2 + p ( 1 + p x_1 ^2 + 2q x_1 x_2 + r x_2 ^2) +
O(|x|^3),$$
$$ g_{22} = - a_{12} x_1 - a_{03} x_2 + r ( 1 + p x_1 ^2 + 2q x_1 x_2 + r x_2 ^2) +
O(|x|^3),$$
$$ g_{12} = g_{21} = - a_{21} x_1 - a_{12} x_2 + q ( 1 + p x_1 ^2 + 2q x_1 x_2 + r x_2 ^2) +
O(|x|^3).$$

Note also that at the point in question, $\det (g) = pr - q^2$.  Running MAPLE, we obtain
that the curvature tensor at this point has 
$$ R_{1212} = -2 (pr - q^2) + \qu \big ( p (a_{12} ^2 - a_{21} a_{03} ) 
+ q (a_{30} a_{03} - a_{21} a_{12} ) + r (a_{21} ^2 - a_{30} a_{12} ) \big )/ (pr - q^2),$$
and hence 
that the (sectional) curvature is
$$ R = -2 + {1\over {4 (pr - q^2)^2 }} \big ( p (a_{12} ^2 - a_{21} a_{03} ) 
+ q (a_{30} a_{03} - a_{21} a_{12} ) + r (a_{21} ^2 - a_{30} a_{12} ) \big ) .$$

If we started instead with a point $L_0$ of $W$, say $L_0 ^3 = A > 0$, then the
curvature at the corresponding point of $W_1$ takes the form 
$$ R = -2 + {A\over {4 (pr - q^2)^2 }} \big ( p (a_{12} ^2 - a_{21} a_{03} ) 
+ q (a_{30} a_{03} - a_{21} a_{12} ) + r (a_{21} ^2 - a_{30} a_{12} ) \big ) .$$

Turning now to the formula for the invariant $S$ on page 167 of Sturmfels, we find that for
the cubic 
$$  F(x_0 ,x_1, x_2) = A x_0^3 - 3 x_0 (p x_1 ^2 + 2q x_1 x_2 + r x_2 ^2)
+ a_{30} x_1 ^3 + 3 a_{21} x_1 ^2 x_2 + 3 a_{12} x_1  x_2 ^2 + a_{03} x_2 ^3 ,$$
the basic invariant $S$ takes the form
$$ S = (pr - q^2)^2 + A \big ( p (a_{12} ^2 - a_{21} a_{03} ) 
+ q (a_{30} a_{03} - a_{21} a_{12} ) + r (a_{21} ^2 - a_{30} a_{12} ) \big ) .$$
From this it follows that
$$ R = - {9\over 4} + {S \over {4 (pr - q^2)^2 }}.$$

An easy check verifies that at the point in question, the Hessian $H$ takes the value
$6^3 A (pr - q^2)$, where $F$ takes the value $A$, and so 
$$ R = - {9\over 4} + {{6^6 S \, F(x_0 , x_1, x_2 )^2}\over {4 H(x_0 , x_1, x_2 )^2}}.$$

For an arbitrary point of $W$, we can always make a change of coordinates  
$$ x_0 \mapsto x_0 + \lambda x_1 + \mu x_2 ,\quad x_1 \mapsto x_1 ,
\quad x_2 \mapsto x_2 , $$ 
so that
the cubic takes the form assumed above.
Since $S$ is invariant under such a transformation in $SL (3, \R )$, as is the value of
the Hessian, the formula remains true in general.\medskip

\noindent\it Examples. \quad\rm  There are two cases above when the curvature is constant
on $W_1$.  The first is when $S=0$, when the curvature is $-{9\over 4}$ everywhere.
Note that $S=0$ if and only if the cubic curve is isomorphic (over $\C$) to the unique
elliptic curve with an automorphism of order three, namely the Fermat cubic $x^3 + y^3 + z^3
=0$, or else some singular specialization of this. Compare this with the calculation from 
Example 2 described in the Introduction; in particular the $-9/4$ occurs because it is $- ({n\over
2})^2$. If $S\ne 0$ and the curvature is constant, then we must have that
$F$ and $H$ are proportional.  This says that every point of the curve is an inflexion point,
and so  the curve reduces to three lines.  Assuming that $W_1$ is non-empty and so 
in particular the 
Hessian is not identically zero, the curve consists of three non-concurrent lines.
An easy calculation then shows that $ {{6^6 S \, F^2}\over {4 H^2}} = {9\over 4} $, 
and so $R$ is identically zero.   Both these cases do occur: if the threefold
admits a birational contraction of a divisor down to a point, then it is an 
easy check from [20] that $S=0$; if $X =\P^1 \times \P ^1 \times \P^1$, then
we are clearly in the second case [18].  

\noindent\it Remarks. \quad\rm   The 
discriminant $\Delta$ of the ternary cubic may be negative, zero or positive;
the former case occurs for instance when $S=0$, and the latter two cases occur for 
 various examples of \CY complete intersections in a
product of three projective spaces.  For \Kahler threefolds
with $b_2 = 3$, the invariant
$S$ may be positive, negative or zero.  
An example with $S<0$ 
is for instance given by taking a line with normal bundle $(-1,-1)$ 
on a quintic hypersurface in
$\P ^4$ and blowing it up, and then blowing up a section of the resulting exceptional surface 
$\P ^1 \times \P^1$.  

The formula in (5.1) makes it possible to check whether a given \Kahler threefold 
with $b_2 =3$ does satisfy the suggested bounds on curvature for the \Kahler moduli,
 and in all the examples calculated by the author the answer has been affirmative.  In 
the \CY cases, all the examples have $S\ge 0$, with strict inequality holding in most cases;
when $S >0$, the curvature
condition has in all cases defined a proper subcone of the index cone. 

In the case of \Kahler threefolds with $b_2 =3$ which 
are complete intersection in a product of three projective spaces, 
the author has made extensive computer calculations,
and there is persuasive numerical evidence here 
both for $S>0$ and the 
\Kahler moduli curvature 
always being semi-negative.  The former property has been checked by the author 
using MATHEMATICA for complete intersections $X$ in 
$\P ^5 \times \P ^5 \times \P ^5$, and the latter property for 
complete intersections in $\P ^3 \times \P ^2 \times \P ^2$.  For three projective spaces of arbitrary dimensions,
there is presumably a combinatorial proof for these assertions, although it seems that the problem is quite hard.
\medskip

\noindent\it Example.\rm\quad  To give the flavour of the calculations 
referred to above, let us consider the complete intersection \CY threefold 
$$\pmatrix{{\P^3}&{\bigl\|}&1&1&2&0\cr{\P^2}&{\bigl\|}&1&1&1&0\cr
{\P^2}&{\bigl\|}&0&0&1&2 \cr}.$$  The configuration is such that the 
\Kahler cone has generators 
corresponding to the three factors, and so the \Kahler cone is 
given in coordinates $(x,y,z)$ by all three coordinates being positive.  Running MAPLE, one discovers that 
$S= 4624$ and the curvature of $\calK_1$ is given by a formula whose denominator is a square and whose numerator 
consists of 13 monomials of degree 6 in $x,y,z$, all of 
whose coefficients are negative.  This property of every  
term being negative is a feature of all the computer calculations performed  
(including various examples of \Kahler 4-folds in the product of three projective spaces that the author has also calculated, using the generalised 
form [23] of (5.1) for quartics).
 \medskip
 
By results of Wall and Jupp [24,9], a simply connected 
compact differentiable 6-manifold with
torsion-free cohomology is determined uniquely by a set of invariants, 
these invariants including the integral cubic
form on  second cohomology given by cup product, the integral linear form on 
 second cohomology given by the first Pontryagin class $p_1$, and the third cohomology,
the first two of these invariants being required to satisfy certain 
simple congruence relations (see [18]
for a good summary).  For the manifold to support a \Kahler structure, 
the Hard Lefschetz Theorem says that 
the Hessian of the cubic form cannot be identically zero [18].  
In the case when $b_2 =3$, the only non-trivial cubic forms excluded 
by this condition are those where the corresponding curve consists of 
three concurrent lines (and 
degenerations of this).  One might think that, even in the case $b_2 =3$, the  
suggested bounds on the 
curvature of $\calK _1$ would rule out some 
further cubic forms --- 
this however is not the case.

\proclaim Proposition 5.2.  For any real cubic form $F(X_0 , X_1 , X_2)$ with Hessian 
$H(X_0 , X_1 , X_2)$
not identically zero, there is a non-empty open subcone of the index cone $W$ for
which the
curvature at the corresponding points of $W_1$ lies between $-3$ and zero.\par

\noindent \it Proof.\quad\rm  The reader is left as an exercise to check the cases when the 
corresponding complex projective curve is reducible.
Suppose now the curve is irreducible (therefore either an elliptic curve, a nodal
cubic, or a cuspidal cubic); it is then well-known that the curve contains at least one
real inflexion point (in the case of a real elliptic curve, there will in fact be three).
Such a real inflexion point is given by a transverse intersection of 
the real curves given by $F$ and $H$.  If we consider the complement of these two real
curves in a small enough neighbourhood of the inflexion point in the affine plane, 
we obtain 
four regions, determined by the signs of $F$ and $H$.  On the region for which $F <0$ and
$H>0$, the signature at the corresponding point of $\R^3$ has to be $(-,-,+)$; we use the
convention here that the first entry indicates the sign of $F$.   As one passes  through
the curve
$F=0$, the signature does not change, only the sign of $F$; we therefore  have locally a
region where the corresponding points of $\R^3$ have signature $(+,-,-)$, that is lie in
the index cone.  If however, we take points in this region away from the inflexion point but
sufficiently close to the curve $F=0$, then the formula in (5.1) ensures that the curvature
condition is also satisfied, and the Proposition follows.  Similarly, points in the region
$F<0$, $H<0$, with
$F$ sufficiently small, will by taking negatives of the vectors in $\R^3$ 
also correspond to points
in $W$ of the type desired.\medskip

 Let us take a specific cubic for illustrative purposes, namely the real nodal
cubic $ F(X_0 , X_1 , X_2 ) = X_1 ^3  + X_0 X_1 ^2 - X_0 X_2 ^2$, or in affine coordinates
$x^3 + x ^2 - y^2$.    By [20], we check that the invariant $S = 1/81$.  The affine form 
of the Hessian is checked easily to be $8 x^2 - 8 y^2 - 24 x y^2$, whose 
zero locus is given by the equation $y^2 = x^2 /(3x+1)$.  The two curves are illustrated 
in the Figure, the solid line denoting the nodal cubic and the broken line its Hessian
curve.

 \midinsert $${\hfil \epsfxsize = 0.35 \hsize \epsfbox{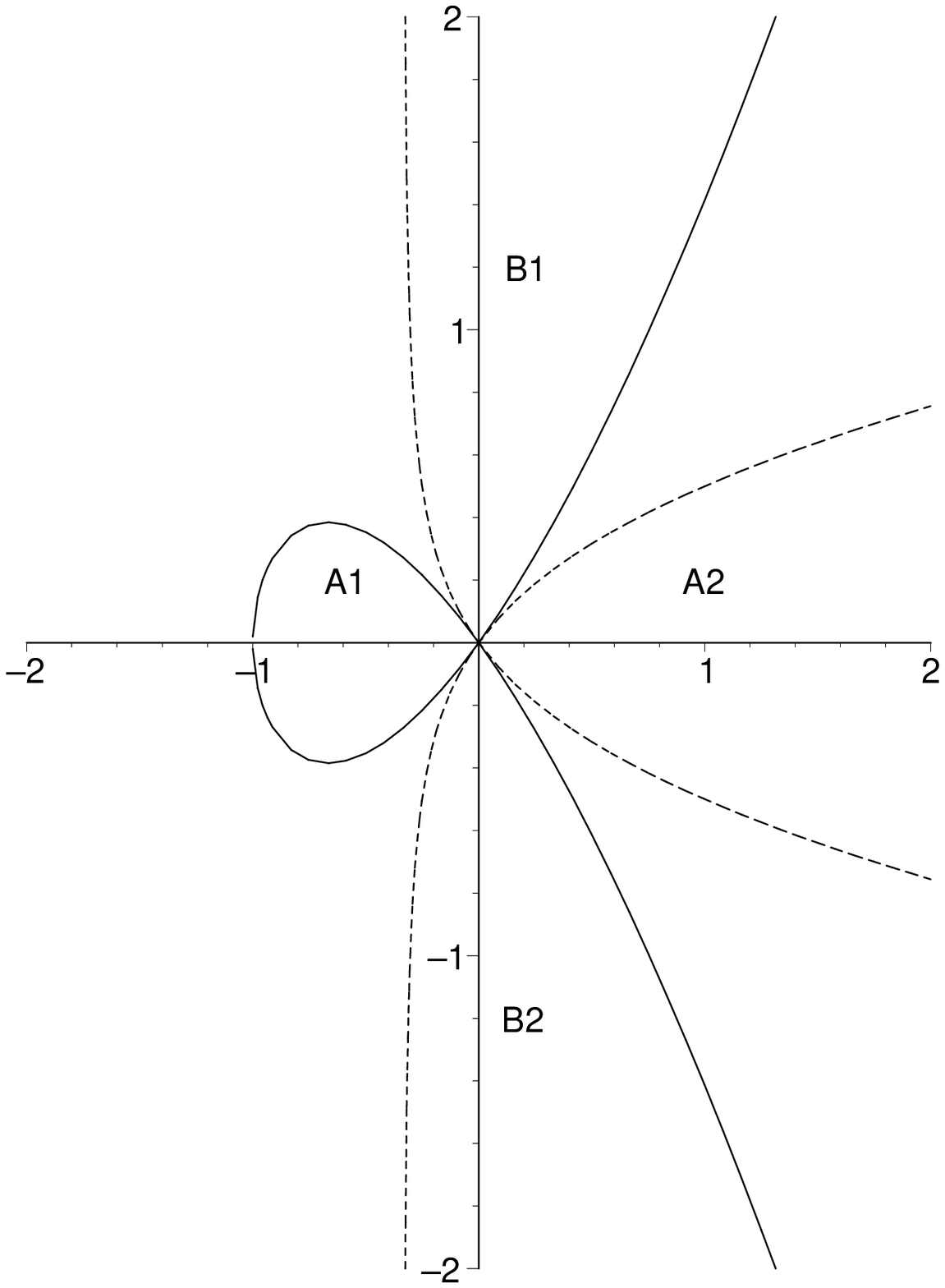} \hfil }$$ \bigskip
 \centerline{The nodal cubic and its Hessian} \medskip \endinsert
We are looking
for  points for which the signature (with the above convention) is either $(+,-,-)$ or
$(-,+,+)$, and these fall in the four regions of the plane shown
(bounded by real branches of the two curves): for regions A1 and
A2 the signature is $(+,-,-)$, and for regions B1 and B2 it is $(-,+,+)$.  

We now
impose the condition that the curvature is non-positive at the corresponding point of
$W_1$; the criterion for this given by (5.1) is invariant on passing 
between a point and its negative
in $\R ^3$.  Thus by (5.1), the curvature condition may be rephrased as the polynomial
$$ (x^3 + x^2 - y^2 )^2 - (x^2 - y^2 - 3 x y^2 )^2 = x \big ( x^5 + 2(x^2 - y^2)(3 y^2 +
x^2) - 9 x y^4 \big )$$
being non-positive.  An elementary calculation shows the real part of the zero locus of this
polynomial consists of  the $y$-axis, together with the curve 
$ y^2 = x^2 (x+2)/(2+3x).$  From this, it is a routine check that the points $(x,y)$ which
satisfy both the index and curvature conditions fall in three regions, whose interiors are as
follows:

\noindent (a)\quad The bounded open set A1 of points with negative $x$-coordinate for which the 
cubic is positive.

\noindent (b)\quad  The two open subsets of B1 and B2 consisting
of points $(x,y)$ in the plane with positive $x$-coordinate and 
for which the 
cubic is negative.  

We should however also take into account the linear form given by $p_1$; recall that for a
complex manifold $X$, we have $p_1  = c_1 ^2 - 2 c_2$, with $c_1, c_2$ denoting the Chern 
classes. 
Let us take the linear form $p_1$  on $\R^3$ to be a negative integral multiple
of $X_1$; note that its zero locus 
corresponds to the
$y$-axis in the affine plane, and therefore separates the open set in (a) from the open sets
in (b).  In the list of invariants given in (1.1) of [18], we may take $w_2 = 0$, $\tau =0$, 
and $b_3$ an arbitrary (even) positive integer.
By taking the cubic form to be a suitable positive integral multiple 
of $X_1 ^3  + X_0 X_1 ^2 - X_0 X_2 ^2$, and similarly with the linear form $p_1$, we may
ensure that the  congruency relations for these invariants to be represented by a 
simply connected smooth 6-manifold with torsion-free homology are
automatically satisfied --- for details of this, see [18].  The smooth manifold $X$ is then
uniquely determined.  We can ask whether it admits any \CY structures.

\proclaim Proposition 5.3.  If the conjectured semi-negativity of the \Kahler moduli curvature 
for \CY threefolds is true,
then the smooth 6-manifolds
$X$ determined by the  above choices of invariants (with $b_3$ arbitrary) carry no \CY 
structures.\par

\noindent \it Proof.\quad\rm  Suppose  $X$ did admit such a \CY structure.
By Lemma 4 of [18], $h^{2,0} =0$ (since otherwise the cubic curve 
consists of either a conic with a transversal line, or three non-concurrent lines),
and so $X$ is projective.
Given that $c_1$ is trivial, the linear form on second cohomology
defined  by $c_2$ is the same as that defined by $- \ha p_1$.  
Since $X$ is minimal, a standard result [17] 
ensures that $c_2$ is non-negative on any \Kahler class, and hence $p_1$ will be
non-positive.  This contradicts our choice of $p_1$, which was 
chosen so as to be strictly positive on every possible
\Kahler class allowed by the conjectured semi-negativity of
the \Kahler moduli curvature.
\medskip
\noindent\it Remarks.\quad\rm  The careful choice of $p_1$ was crucial in the above argument, 
since the nodal cubic does occur for the Calabi--Yau threefold 
$$\pmatrix{{\P^2}&{\bigl\|}&1&2\cr{\P^2}&{\bigl\|}&2&1\cr
{\P^1}&{\bigl\|}&0&2 \cr},$$ whose cubic form, a multiple of $4 x^2 y + 2 x y^2
 + 2x^2 z  + 
2 y^2 z + 10 x y z$, is irreducible and has a node at $(0:0:1)$.  This incidentally answers 
a question raised in Section 5.3 of [18].  Additionally, the 
semi-negativity of the \Kahler moduli curvature is known in this case.\par 
We chose the nodal cubic for the above example because 
the equations turn out to be relatively simple.  A slightly more complicated example
would be to take the cubic form $x^3 - x +1 - y^2$, with $S= 1/27$, corresponding to a
smooth  elliptic curve.  Again, one finds that there are three open subsets in the plane
that could correspond to points of $\calK _1$ where the curvature was non-positive.  The
(rational) line through the three real  inflexion points of the elliptic curve 
separates one of these subsets (on which both the cubic and its
Hessian are positive) from the other two (on which both the cubic and its
Hessian are negative).  The same argument as above then goes through. 
\bigskip

\noindent\bf References.\hfill\rm

\noindent 1.  Besse, A.L.:  Einstein Manifolds.  Berlin-Heidelberg: Springer 1987.

\noindent 2.  Candelas, P., de la Ossa, X.C., Green, P.S., Parkes, L.:  A pair of \CY manifolds as an exactly soluble superconformal theory.  In : 
Yau, S.-T. (ed.) : Essays in mirror manifolds, pp. 31-95.  
International Press 1992.

\noindent 3.  Fang, H., Lu, Z.:  Generalized Hodge metrics and BCOV 
torsion on \CY moduli.  Preprint, 2003. 
ArXiv:math. DG/0310007 (2003).

\noindent 4.  Griffiths, P., Harris, J.:  Principles of Algebraic Geometry.
New York : Wiley 1978.

\noindent 5.  Gross, Mark:  Special Lagrangian Fibrations II: Geometry.  A survey of
Techniques in the study of Special Lagrangian Fibrations.  In : Yau, S.-T. (ed.)
: Surveys in Differential Geometry V, pp. 341-403. 
International Press, Somerville MA: International Press 1999.

\noindent 6.  Hitchin, N.:  The moduli space of special Lagrangian submanifolds.
Ann. Scuola Norm. Sup. Pisa Cl. Sci. \bf 25\rm , 503-515 (1997).

\noindent 
7.  Huybrechts, D.: Products of harmonic forms and rational curves.  Doc. Math. \bf 6\rm
, 227-239 (2001).

\noindent
8.  Jost, J.:  Riemannian Geometry and Geometric Analysis.  
Berlin-Heidelberg: Springer 2002.

\noindent 9.  Jupp, P.: Classification of certain 6-manifolds.  Math. Proc. Cam. Phil. Soc.
\bf 73\rm , 293-300 (1973).

\noindent 10.  Kobayashi, S.:  Transformation Groups in Differential Geometry.  
Berlin-Heidelberg: Springer 1985.

\noindent 11.   Kobayashi, S., Nomizu, K.:  Foundations of Differential Geometry.
New York: Wiley 1963
 
\noindent 12.  Leung, N.C.:  Mirror Symmetry without corrections.  Comm. Anal. Geom., to appear.
ArXiv:math. DG/0009235.

\noindent 13.  Liu, K., Todorov, A., Yau, S.-T., Zuo, K.:
Shafarevich's Conjecture on CY manifolds I (Moduli of CY 
manifolds).  Preprint, 2003. ArXiv:math. AG/0308209.

\noindent
14.  Lu, Z.:  On the Hodge metric of the universal deformation space of \CY 
threefolds.  J. Geom. Anal. \bf 11\rm , 103-118 (2001).

\noindent
15.  Lu, Z., Sun, X.:  Weil--Petersson geometry on the moduli space of polarized \CY manifolds.  
Preprint 2002 (to appear in Journal de l'Institut Mathematique de Jussieu).

\noindent
16.  Miller, T.J.:  On the formality of $(k-1)$ connected compact manifolds of dimension 
less than or equal to $(4k-2)$.  Illinois J. Math. \bf 23\rm , 253-258 (1979).

\noindent
17. Miyaoka, Y.:  The Chern classes and Kodaira dimension on a minimal variety.  In :
Oda, T. (ed.) Algebraic Geometry, Sendai 1985.  Adv. Stud. Pure Math. Vol. 10, pp. 
449-476.  Amsterdam: North Holland 1987.

\noindent
18.  Okonek, Ch., Van der Ven, A.:  Cubic forms and complex 3-folds.  Enseign. Math. (2) 
\bf 41\rm , 297-333  (1995).

\noindent 19.  Sakai, T.:  Riemannian Geometry.  Providence: AMS 1996.

\noindent 20. Sturmfels, B.:  Algorithms in Invariant Theory.  New York: Springer 1993.

\noindent 21.  Tian, G.:  Smoothness of the universal deformation space of compact \CY 
manifolds and its Weil--Petersson metric.  In : Yau, S.-T. (ed.) Mathematical Aspects of
String Theory, pp. 629-646.  
Singapore: World Scientific 1981.

\noindent 22.  Todorov, A.N.: The Weil--Petersson geometry of the moduli space of 
$SU(n\ge 3)$ (Calabi--Yau) manifolds I.  Commun. Math. Phys. \bf 126\rm , 325-346
 (1989).

\noindent 23.  Totaro, B.:  The curvature of a Hessian metric.  Int. J. Math., to appear.  ArXiv:math. DG/0401381.

\noindent 24.  Wall, C.T.C.:  Classification Problems in Differential
Topology V.  On  certain 6-Manifolds.  Invent. math. \bf 1\rm , 139-155
(1966).

\noindent 25.  Wang, C.L.:  On the incompleteness of the Weil--Petersson
metric along degenerations of \CY manifolds.  Mathematical Research
Letters \bf 4\rm , 157-171 (1997).

\noindent
26.  Wang, C.L.:  Curvature properties of the \CY moduli.  Doc. Math. \bf 8\rm , 577-590 (2003).

\noindent 27.  Wilson, P.M.H.:  The \Kahler cone on \CY threefolds. 
Invent. math. \bf 107\rm , 561-583 (1992).  Erratum : Invent. math. \bf
114\rm
 , 231-233 (1993).

\noindent 28.  Zheng, F.:  Complex Differential Geometry.
AMS/IP Studies in Adv. Math. \bf 18. \rm 
 Providence: AMS/IP 2000.
\vfill

\end